\documentclass[12pt]{article}
\usepackage[latin1]{inputenc}
\usepackage{amssymb,amsmath,amsfonts}
\usepackage{xcolor}
\newtheorem{theorem}{Theorem}
\newtheorem{proposition}[theorem]{Proposition}
\newtheorem{lemma}[theorem]{Lemma}

\newtheorem{fact}[theorem]{Fact}

\newcommand{\R}{\mathbb{R}}
\newcommand{\U}{\mathcal{U}}

\newcommand{\E}{{\cal E}}
\newcommand{\spa}{\mbox{span\,}}
\newcommand{\Ima}{\mbox{Im\,}}
\newcommand{\rank}{\mbox{rank }}
\newcommand{\kerl}{\mbox{ker}}

\newcommand{\po}{{\hspace*{-1ex}}{\bf .  }}

\newcommand{\nap}{\nabla^{\perp}}
\newcommand{\nab}{\tilde\nabla}
\def\lp{{\langle\!\langle}}\vspace{2ex}
\def\rp{{\rangle\!\rangle}}
\def\<{{\langle}}
\def\>{{\rangle}}
\def\Sal{{\cal S}}
\def\J{{\cal J}}
\def\T{{\cal T}}

\def\a{\alpha}
\def\be{\begin{equation} }
\def\ee{\end{equation} }

\def\proof{\noindent\emph{Proof: }}
\def\qed{\ifhmode\unskip\nobreak\fi\ifmmode\ifinner
\else\hskip5 pt \fi\fi\hbox{\hskip5 pt \vrule width4 pt
height6 pt  depth1.5 pt \hskip 1pt }}
\newcommand\blfootnote[1]{
\begingroup
\renewcommand\thefootnote{}\footnote{#1}
\addtocounter{footnote}{-1}
\endgroup
}
\begin{document}

\title{Non-holomorphic Kaehler submanifolds\\ of Euclidean space}
\author{Sergio Chion and Marcos Dajczer}
\date{}
\maketitle

\begin{abstract} This paper is about non-holomorphic 
isometric immersions of Kaehler manifolds into Euclidean 
space $f\colon M^{2n}\to\R^{2n+p}$, $p\leq n-1$, with low 
codimension $p\leq 11$. In particular, it addresses a 
conjecture proposed by J. Yan and F. Zheng. The claim  
that if the index of complex relative nullity of the 
submanifold satisfies $\nu_f^c<2n-2p$ at any point,
then $f(M)$ can be realized as a holomorphic submanifold 
of a  non-holomorphic Kaehler submanifold of $\R^{2n+p}$ 
of larger dimension. This conjecture had previously been 
confirmed 
by  Dajczer-Gromoll for codimension $p=3$, and then by 
Yan-Zheng for $p=4$. For codimension $p\leq 11$,
we already showed that the pointwise structure of the 
second fundamental form of the submanifold aligns 
with the anticipated characteristics, assuming the 
validity of the conjecture. In this paper, we confirm 
the conjecture until codimension $p=6$, whereas for 
codimensions $7\leq p\leq 9$ it is also possible that 
the submanifold exhibits a complex ruled structure 
with rulings of a specific dimension.  Moreover, we 
prove that the claim of the conjecture holds for 
codimensions $7\leq p\leq 11$ albeit subject to an 
additional assumption.
\end{abstract}
\blfootnote{\textup{2020} \textit{Mathematics Subject Classification}:
53B25, 53C42.}

\section{Introduction}

In a seminal work, E. Calabi \cite{Ca} characterized the real 
analytic Kaehler metrics that can be locally isometrically and 
holomorphically immersed into complex Euclidean space. 
This is done by means of a metric condition which only a highly 
restricted subset of such metrics meet. In addition, it was  proved 
by Calabi that these immersed submanifolds are isometrically rigid 
regardless of their codimension.

As for the Kaehler manifolds that admit isometric immersions 
with low codimension into Euclidean space, it has been shown
that the submanifolds are generically holomorphic. A first 
result in that direction was achieved in \cite{DR}, focusing 
on submanifolds laying in any codimension, provided it does 
not exceed one third of the dimension of the manifold. A stronger 
result was given in \cite{CCD} for submanifold of lower dimensions 
concerning codimension. However, its applicability is constrained 
to situations where the latter does not surpass eleven.

Extending beyond the confines of holomorphicity, isometric 
immersions of Kaehler manifolds into Euclidean space emerged 
as a natural focus for further investigation following the 
developments in \cite{DG1}. In broad terms, this paper shows 
that non-holomorphic minimal Kaehler submanifolds are the 
ones permitting isometric deformations while preserving the 
Gauss map. Note that the outcomes in \cite{DG1} or \cite{DT} 
established the persistence of the metric condition specified 
by Calabi for immersions that are simply minimal. This is due 
to the fact that any such submanifold can be realized as the 
real part of its holomorphic isometric representative. 

The result in \cite{CCD} shows that the Kaehlerness of the 
manifold imposes severe restrictions on a non-holomorphic 
submanifold. This was already well-known for very low 
codimension and is further solidified by the results 
in our current paper. The main finding in this paper is, 
in essence, that in codimension of at most $9$ any 
non-holomorphic Kaehler submanifold must be either complex 
ruled or a holomorphic submanifold of a complex ruled 
submanifold of higher dimension. The precise result 
substantiating this rather succinct statement 
is elaborated upon in the subsequent discussion.
\vspace{1ex}

An isometric immersion $f\colon M^{2n}\to\R^{2n+p}$ is called 
a \emph{real Kaehler submanifold} if $(M^{2n},J)$ is a nowhere 
flat connected Kaehler manifold with complex 
dimension $n\geq 2$ immersed into standard flat Euclidean space 
with substantial codimension. The latter condition means 
that the codimension of the submanifold cannot be reduced 
even locally. In other words, no open subset of the submanifold 
is contained in an affine proper subspace of the Euclidean ambient 
space. Moreover, if $p$ is even, we require that $f$ restricted 
to any open subset of $M^{2n}$ not be holomorphic with respect 
to any possible complex structure of $\R^{2n+p}$. 
An abundance of basic information on the subject of real Kaehler 
submanifolds is provided in Chapter $15$ in \cite{DT}. This 
encompasses a parametric local classification for the real Kaehler 
hypersurfaces and, in particular, a Weierstrass representation
in the minimal case. In the case where minimality is not assumed, 
these submanifolds may lack the property of being real analytic.
For submanifolds in codimension two we refer to the Appendix in 
\cite{CD1} where, in particular, other classifications results 
are discussed.

In the pursuit of classifying real Kaehler submanifolds that 
lie in codimension greater than two, a crucial preliminary step 
involves excluding the following compositions of isometric 
immersions: Let $F\colon N^{2n+2m}\to\R^{2n+p}$, $m\geq 1$, 
be a real Kaehler submanifold and let 
$j\colon M^{2n}\to N^{2n+2m}$ be an holomorphic 
submanifold. Then, the composition  of isometric immersions 
$f=F\circ j\colon M^{2n}\to\R^{2n+p}$ is, generically, a real 
Kaehler submanifold. Concomitantly, we say that a  real 
Kaehler submanifold $f\colon M^{2n}\to\R^{2n+p}$ has a 
\emph{Kaehler extension} if there is a real 
Kaehler submanifold $F\colon N^{2n+2m}\to\R^{2n+p}$, $m\geq 1$, 
and a holomorphic isometric embedding $j\colon M^{2n}\to N^{2n+2m}$ 
such that $f=F\circ j$.

A local characterization of the real Kaehler submanifolds that admit 
a Kaehler extension have been achieved for codimension $p=3$ by 
Dajczer and Gromoll \cite{DG2}, and then for codimension $p=4$ by Yan 
and Zheng \cite{YZ}. For the latter case, we provided in \cite{CD1} 
some additional information. These results have been obtained 
under an assumption on the dimension of the complex kernel of the 
second fundamental form of the submanifold discussed next. 
\vspace{1ex}

Let $\a\colon TM\times TM\to N_f M$ denote the second fundamental 
form of a real Kaehler submanifold $f\colon M^{2n}\to\R^{2n+p}$. 
The \emph{index of complex relative nullity} of $f$ at $x\in M^{2n}$ 
is $\nu_f^c(x)=\dim\mathcal{N}(\a)(x)\cap J\mathcal{N}(\a)(x)$, 
where
$$ 
\mathcal{N}(\a)(x)
=\{Y\in T_xM\colon \a(X,Y)=0\;\;\mbox{for any}\;\;X\in T_xM\}
$$
is the relative nullity subspace of $f$ at $x\in M^{2n}$.
\vspace{1ex}

A conjecture formulated by Yan and Zheng \cite{YZ} applies to 
real Kaehler submanifolds $f\colon M^{2n}\to\R^{2n+p}$, $p\leq n-1$, 
of codimension of at most $p=11$ that satisfy the condition 
$\nu_f^c(x)<2n-2p$ at any point $x\in M^{2n}$. Their claim is 
that $f$ will necessarily admit a Kaehler extension.
The assertion of this bold conjecture aligns with what has 
already been established for codimensions $3$ and $4$ 
in the aforementioned papers. 
It is noteworthy to observe that an affirmative resolution 
to the conjecture holds significant implications for the metric 
of the manifold. Specifically, in cases where the submanifold 
does not allow such extension, the curvature tensor would exhibit 
a complex kernel with a minimum dimension of $\nu_f^c$ at any point.

In a prior work we took a significant step towards addressing 
the conjecture by Yan and Zheng. Specifically, we showed in 
\cite{CD2} that the pointwise behavior of the second fundamental 
form of the submanifold aligns with the expectations if the 
conjecture holds true. This result is given below as 
Theorem \ref{alglemma} and is of fundamental use in this paper.
\vspace{1ex}

A real Kaehler submanifold $f\colon M^{2n}\to\R^{2n+p}$ is said 
to be $R^\rho$-\emph{complex ruled} (or just $\rho$-complex 
ruled if convenient) if there is a holomorphic foliation 
$R^\rho$ of $M^{2n}$ with leaves of real even dimension $\rho$ such 
that the image  by $f$ of any leaf is part of an affine vector 
subspace of $\R^{2n+p}$. In other words, the submanifold  $f(M)$ 
is a holomorphic bundle of affine vector subspaces. Notice that
a $\rho$-complex ruled submanifold may actually carry a complex 
ruled foliation of leaves of higher dimension. Observe also that
if $\nu_f^c(x)=m>0$ at any $x\in M^{2n}$ then $f$ is $m$-complex 
ruled by the integral leaves of the relative nullity distribution.
\vspace{1ex}

In this paper, we establish the validity of the Yan-Zheng conjecture 
up to codimension $p=6$, whereas for codimensions $7\leq p\leq 9$ we 
argue that a second possibility may occur, namely, that the submanifold 
can be complex ruled. As for codimensions $p=10$ and $p=11$, 
we affirm the conjecture's veracity, albeit under an additional 
assumption.  
In the concluding part of this introduction, we discuss why in 
codimensions higher than $6$ the emergence of complex ruled 
submanifolds is for us an expected alternative.
\vspace{1ex}

The following is the main result of this paper.

\begin{theorem}\po\label{th1}
Let $f\colon M^{2n}\to\R^{2n+p}$, $p\leq n-1$ and $3\leq p\leq 9$, 
be a real Kaehler submanifold with $\nu_f^c(x)<2n-2p$ at any 
$x\in M^{2n}$. Assume for $7\leq p\leq 9$ that $f$ restricted 
to any open subset of $M^{2n}$ is not $(2n-t)$-complex ruled 
for $(p,t)=(7,8),(8,12)$ or $(9,16)$. Then $f$ restricted 
to any connected component of an open dense subset of $M^{2n}$ 
has a Kaehler extension $F\colon N^{2n+2r}\to\R^{2n+p}$ where 
$2r<p$ and such that \mbox{$\nu_F^c(z)\geq 2n-2p+6r$} at any 
$z\in N^{2n+2r}$.  Moreover, the extension $F$ is minimal if 
and only if $f$ is minimal.
\end{theorem}

Observe that the above result includes the class of minimal 
submanifolds, a situation not covered by the Yan-Zheng 
conjecture. The Kaehler extension in the above statement, 
as well as all the others in this paper, is explicitly 
constructed. Unfortunately, the size of our estimates 
impedes the extension of the above result to codimensions 
$10$ and $11$, but that remains a viable possibility.

In \cite{CG} the authors assert the verification of the 
Yan-Zheng conjecture for codimensions $5$ and $6$. Their 
proof is intricately linked to Proposition $2$ in the paper. 
Regrettably, this result is flawed, as the provided ruled 
extension may not necessarily be a Kaehler manifold. 
Indeed, for their extension the Kaehler property holds 
true exclusively in a very specific case, coinciding with 
our extension in this particular instance.
\vspace{1ex}

Our next result shows that a real Kaehler submanifold with
codimension $7\leq p\leq 11$ necessarily has to admit a Kaehler 
extension upon an additional assumption concerning the structure 
of the second fundamental form.
\vspace{1ex}

The \emph{index of complex $s$-nullity} $\nu^c_s(x)$, 
$1\leq s\leq p$, at $x\in M^{2n}$ of a real Kaehler submanifold 
$f\colon M^{2n}\to\R^{2n+p}$ is defined by 
$\nu^c_s(x)=\max_{U^s\subset N_fM(x)}\nu^c(\a_{U^s})$,
where $U^s$ is any $s$-dimensional vector subspace,  
$\a_{U^s}=\pi_{U^s}\circ\a$ being $\pi_{U^s}\colon N_fM\to U^s$ 
the projection and 
$\nu^c(\a_{U^s})(x)=\dim\mathcal{N}_c(\a_{U^s}(x))$ with
$$
\mathcal{N}_c(\a_{U^s}(x))
=\mathcal{N}(\a_{U^s}(x))\cap J\mathcal{N}(\a_{U^s}(x))
$$
and
$$
\mathcal{N}(\a_{U^s}(x))
=\{Y\in T_xM\colon \a_{U^s}(X,Y)=0\;\mbox{for any}\;X\in T_xM\}.
$$
Notice that $\nu^c_f(x)=\nu_p^c(\a)(x)$. Moreover, that we have
$\nu^c_s(x)\geq\nu^c_{s+1}(x)$ for $1\leq s\leq p-1$.
\vspace{1ex}

In our prior research \cite{CDa}, we established that real 
Kaehler submanifolds, irrespective of their codimension, 
must necessarily be minimal if specific upper bounds on 
the indexes of complex $s$-nullities are satisfied. However, 
after relaxing some of these assumptions, a distinct result 
emerges for codimensions within the range of $7$ to $11$.

\begin{theorem}\po\label{th2}
Let $f\colon M^{2n}\to\R^{2n+p}$, $p\leq n-1$ and $7\leq p\leq 11$, 
be a real Kaehler submanifold satisfying $\nu_s^c(x)<2n-2s$ 
at any $x\in M^{2n}$, where 
$$
(p,s)=(7,5;7),\, (8,6;8),\,(9,5;7;9),\,(10,6;8;10)\;\text{and}
\;(11,5;7;9;11).
$$
Then $f$ restricted to any connected component of an open 
dense subset of $M^{2n}$ has a Kaehler 
extension $F\colon N^{2n+2r}\to\R^{2n+p}$ which satisfies that 
$\nu^c_F(z)\geq 2n-2p+6r$ at any $z\in N^{2n+2r}$, where 
$(p,r)=(7,2\,\text{or }3),\,(8,2\,\text{or }3),\,(9,3\,\text{or }4),
\,(10,3\,\text{or }4)\;\text{or}\;(11,4\,\text{or }5)$. 
Moreover, the extension $F$ is minimal if and only if $f$ is 
minimal.
\end{theorem}

The technique employed  in this paper is strongly inspired 
by the one introduced in \cite{DF} to treat the genuine 
rigidity of submanifolds in space forms. Further refinement 
and enhancement of this technique is used in Chapter $12$ 
of \cite{DT}. 

An isometric deformation of a submanifold in Euclidean space is 
said to be genuine if it cannot be obtained induced by the 
deformation of a submanifold of larger dimension in which it is 
isometrically embedded. It was first observed in \cite{CG} that 
there is a relation between the subject of this paper and genuine 
rigidity. This is connected to the fact that any simply connected 
minimal Euclidean Kaehler submanifold admits an associated one 
parameter family of isometric deformations of the same type.  
That family is trivial if and only if the submanifold is 
holomorphic. In the study of submanifolds that are genuinely 
rigid the possible presence of ruled examples also shows up; see 
Theorem $12.10$ in \cite{DT}. We encounter a somewhat analogous 
scenario here since, in both instances, the actual existence of 
these special ruled submanifolds remains a challenging open problem.
In the present case, the construction of such an example
in codimension $7$ to $11$ would contradict the conjecture
posited by Yan and Zheng.

\section{Preliminaries}

In this section, we first introduce some notations and present 
several facts to be referenced throughout the paper. Then we 
recall the outcome from \cite{CD2} utilized to elucidate the 
pointwise structure of the second fundamental form of a real 
Kaehler submanifold in codimension up to eleven.
\vspace{2ex}

Let $\varphi\colon V_1\times V_2\to W$ be a bilinear form 
between finite dimensional real vector spaces. The image 
$\mathcal{S}(\varphi)$ of $\varphi$ is the vector subspace 
of $W$ defined by
$$
\mathcal{S}(\varphi)=\spa\{\varphi(X,Y)
\;\mbox{for all}\;X\in V_1\;\mbox{and}\;Y\in V_2\}.
$$
The (right) nullity $\mathcal{N}(\varphi)$ of $\varphi$ 
is the vector subspace of $V_2$ given by
$$
\mathcal{N}(\varphi)=\{Y\in V_2\colon\varphi(X,Y)=0
\;\mbox{for all}\;X\in V_1\}
$$
whose dimension $\nu(\varphi)$ is called the \emph{index of 
nullity} of $\varphi$.
\vspace{1ex}

A vector $X\in V_1$ is called a (left) \emph{regular element} 
of $\varphi$ if $\dim\varphi_X(V_2)=\kappa(\varphi)$, where
$$
\kappa(\varphi)=\max_{X\in V_1}\{\dim\varphi_X(V_2)\}
$$
and $\varphi_X\colon V_2\to W$ is the linear map defined by 
$\varphi_XY=\varphi(X,Y)$.  The subset of $V_1$ of regular 
elements of $\varphi$ is denoted by $RE(\varphi)$. 

\begin{proposition}\label{facts}\po The following facts hold:
\begin{itemize}
\item[(i)] The subset $RE(\varphi)\subset V_1$ is open 
and dense.
\item[(ii)] Let $W$ be endowed with an inner product. Set 
$\tau_\varphi(X)=\dim\varphi_X(V_2)\cap\varphi_X(V_2)^\perp$
and $\tau(\varphi)=\min_{X\in RE(\varphi)}\{\tau_\varphi(X)\}$.
Then
$$
RE^\#(\varphi)=\{X\in RE(\varphi):\tau_\varphi(X)=\tau(\varphi)\}
$$ 
is an open dense subset of $V_1$. 
\end{itemize}
\end{proposition}

\proof Part $(i)$ is the Proposition $4.4$ in \cite{DT}.  
The proof of Lemma $2.1$ in \cite{DR} or the  Proposition 
$4.18$ in \cite{DT} give part $(ii)$.\vspace{2ex}\qed

The bilinear form $\varphi$ is said to be  \emph{flat} with
respect to an inner product $\<\,,\,\>$ on $W$ if 
$$
\<\varphi(X,Y),\varphi(Z,T)\>-\<\varphi(X,T),\varphi(Z,Y)\>=0
$$
for any $X,Z\in V_1$ and $Y,T\in V_2$. 

\begin{proposition}\po If $X\in RE(\varphi)$ and $N(X)
=\ker\varphi_X$ then
\be\label{firstst}
\Sal(\varphi|_{V_1\times N(X)})\subset\varphi_X(V_2).
\ee
If $\varphi$ be flat with respect to an inner product on $W$ 
then
\be\label{secondst}
\Sal(\varphi|_{V_1\times N(X)})\subset\varphi_X(V_2)
\cap\varphi_X(V_2)^\perp.
\ee
\end{proposition}

\proof This is Proposition $4.6$ in \cite{DT}.\vspace{2ex}\qed

In the sequel $U^p$ stands for a $p$-dimensional vector space 
endowed a positive definite inner product $\<\,,\,\>$.  
Let $W^{p,p}=U^p\oplus U^p$ and $\pi_1\colon W^{p,p}\to U^p$ 
(respectively, $\pi_2$) denote taking the first (respectively, 
second) component of $W^{p,p}$. 
Let $W^{p,p}$ be endowed with the inner product 
$\lp\,,\,\rp$ of signature $(p,p)$ given by
$$
\lp(\xi_1,\xi_2),(\eta_1,\eta_2)\rp
=\<\xi_1,\eta_1\>-\<\xi_2,\eta_2\>.
$$
Then $W^{p,p}$ carries the complex structure 
$\mathcal{T}\in\text{Aut}(W)$ defined by 
\be\label{mathcalT}
\mathcal{T}(\xi,\eta)=(\eta,-\xi)
\ee
which satisfies that
$\lp\mathcal{T}\delta,\nu\rp=\lp\delta,\mathcal{T}\nu\rp$.
\vspace{1ex}

A vector subspace $L\subset W^{p,p}$ called \emph{degenerate} 
if $L\cap L^\perp\neq 0$ and  \emph{nondegenerate} if otherwise.  
A degenerate vector subspace $L\subset W^{p,p}$ is called 
\emph{isotropic} if $L=L\cap L^\perp$. 
\vspace{1ex}

Now let $V_2$ be endowed with a complex structure, that is, 
there is $J\in\text{Aut}(V_2)$ such that $J^2=-I$. 
Assume that the bilinear form 
$\varphi\colon V_1\times V_2\to W^{p,p}$ satisfies 
$$
\mathcal{T}\varphi(X,Y)=\varphi(X,JY)\;\,\mbox{for any}\;\,
X\in V_1\;\mbox{and}\;Y\in V_2.
$$
The isotropic vector subspace 
$\U=\Sal(\varphi)\cap\Sal(\varphi)^\perp$ 
satisfies $\mathcal{T}\U=\U$. In fact, if we have that
$\lp\varphi(X,Y),(\xi,\bar\xi)\rp=0$ for any $X\in V_1$ 
and $Y\in V_2$, then 
$$
\lp\varphi(X,Y),\mathcal{T}(\xi,\bar\xi)\rp
=\lp\mathcal{T}\varphi(X,Y),(\xi,\bar\xi)\rp
=\lp\varphi(X,JY),(\xi,\bar\xi)\rp=0.
$$

\begin{proposition}\po\label{even}
The following facts hold:
\begin{itemize}
\item[(i)] $\mathcal{T}|_{\Sal(\varphi)}\in\text{Aut}\,\,(\Sal(\varphi))$
and $\mathcal{T}|_{\U}\in\text{Aut}(\U)$ are complex structures.
\item[(ii)] The vector subspaces $\Sal(\varphi)$ and $\U$ of 
$W^{p,p}$ have even dimension.
\item[(iii)] The vector subspace $\mathcal{N}(\varphi)\subset V_2$ 
is $J$-invariant and thus of even dimension.
\item[(iv)] If $\Omega=\pi_1(\U)$ then $\dim\Omega=\dim\U$. 
If $\varphi_\Omega=\pi_{\Omega\times\Omega}\circ\varphi$
then $\Sal(\varphi_\Omega)=\U$.
\end{itemize}
\end{proposition}

\proof The proofs of parts $(i)$ to $(iii)$ are immediate. 
Since $\U$ is an isotropic vector subspace then 
$\pi_1|_{\U}\colon\U\to\Omega$ is an isomorphism.  By the 
second statement of part $(i)$ we have $\pi_2(\U)=\Omega$
and thus
$\U\subset\Sal(\varphi_\Omega)\subset\Omega\oplus\Omega$.  
Since $\U$ is isotropic and $\dim\U=\dim\Omega$ then
$\U=\Sal(\varphi_\Omega)$.\vspace{1ex}\qed

Let $\a\colon V^{2n}\times V^{2n}\to U^p$ be a symmetric 
bilinear form and $J\in\text{Aut}(V)$ a complex structure. 
Let $\gamma\colon V^{2n}\times V^{2n}\to W^{p,p}$ be the 
bilinear form defined by 
\be\label{gamma}
\gamma(X,Y)=(\a(X,Y),\a(X,JY)).
\ee
Then let 
$\beta\colon V^{2n}\times V^{2n}\to W^{p,p}$ be the 
bilinear form defined by
\be\label{beta}
\beta(X,Y)=\gamma(X,Y)+\gamma(JX,JY).
\ee

\begin{theorem}\po\label{alglemma}
Assume that the bilinear forms 
$\gamma,\beta\colon V^{2n}\times V^{2n}\to W^{p,p}$,
$p\leq n$, are flat and satisfy
\be\label{productflat2}
\lp\beta(X,Y),\gamma(Z,T)\rp=\lp\beta(X,T),\gamma(Z,Y)\rp
\;\,\mbox{for any}\;\, X,Y,Z,T\in V^{2n}.
\ee
Denote $\Omega=\pi_1(\Sal(\gamma)\cap\Sal(\gamma)^\perp)$ 
and let $U^p=\Omega\oplus P$ be an orthogonal decomposition. 
If $p\leq 11$ and $\nu(\gamma)<2n-\dim\Sal(\gamma)$ 
then $\dim\Omega\geq 2$ and the following facts hold:
\begin{itemize}
\item[(i)] There is an isometric complex structure 
$\J\in\text{End}(\Omega)$ so that $\a_\Omega=\pi_\Omega\circ\a$ 
satisfies 
$$
\J\a_\Omega(X,Y)=\a_\Omega(X,JY)\;\,\mbox{for any}\;\,X,Y\in V^{2n}.
$$  
\item[(ii)] The  bilinear form $\gamma_P=\pi_{P\times P}\circ\gamma$ 
is flat, the vector subspace $\Sal(\gamma_{P})$ is nondegenerate and  
$\nu(\gamma_P)\geq 2n-\dim\Sal(\gamma_P)$.
\end{itemize}
\end{theorem}

\proof This follows from Theorem $9$ in \cite{CD2}.\qed

\section{Complex ruled Kaehler submanifolds}

In this section, we first recall a result from \cite{CD1}. 
There it was established that the presence of a complex 
normal vector subbundle on a real Kaehler submanifold meeting 
certain conditions allows the explicit construction of a 
Kaehler extension.  Subsequently, we analyse the situation 
when such a Kaehler extension turns out to be trivial. 
This leads us to the conclusion that the submanifold has to 
be complex-ruled.

\subsection{Kaehler extension}

Let $f\colon M^{2n}\to\R^{2n+p}$ be a real Kaehler submanifold
with second fundamental form $\a\colon TM\times TM\to N_fM$. 
We make use of the notation $N_1^f(x)=\Sal(\a(x))$ since this 
normal vector subspace is usually called the first normal space
of $f$ at $x\in M^{2n}$. 

Let $L^{2\ell}$ stand for a proper vector subbundle of the normal 
bundle of (real) rank $2\ell>0$  such that $L(x)\subset N_1^f(x)$ 
at any $x\in M^{2n}$. Here and elsewhere we call rank of a vector 
bundle the dimension of its fibers.
Then $L^{2\ell}$ is taken endowed with the induced metric and 
vector bundle connection. We assume that $L^{2\ell}$ 
carries a complex structure $\J\in\Gamma(\text{Aut}(L))$, that is, 
a vector bundle isometry satisfying $\J^2=-I$. We also 
assume that the $J$-invariant tangent vector subspaces 
$D(x)=\mathcal{N}_c(\a_{L^\perp}(x))$ have constant (even) dimension 
$d>0$, and thus form a holomorphic tangent subbundle $D$ of rank $d$. 
In addition, we require the following two conditions:
\begin{itemize}
\item[($\mathcal{C}_1$)] The complex structure $\J\in\Gamma(\text{Aut}(L))$ 
is parallel, that is, 
$$
(\nap_X\J\eta)_L=\J(\nap_X\eta)_L \;\,\mbox{for any}\;\,X\in\mathfrak{X}(M)
\;\mbox{and}\;\eta\in\Gamma(L)
$$ 
and the second fundamental form of $f$ satisfies
\be\label{cond}
\J\a_L(X,Y)=\a_L(X,JY)\;\,\mbox{for any}\;\, X,Y\in\mathfrak{X}(M)
\ee
or, equivalently, that $A_{\J\eta}= J\circ A_\eta
=-A_\eta\circ J\;\,\mbox{for any}\;\,\eta\in\Gamma(L)$.
\item[($\mathcal{C}_2$)] The subbundle $L^{2\ell}$ is parallel 
along $D$ in the normal connection of $f$, that is,
$$
\nap_Y\eta\in\Gamma(L)\;\,\mbox{for any}\;\,Y\in\Gamma(D)\;\mbox{and}\;  
\eta\in\Gamma(L).
$$
\end{itemize}

Let the vector bundle $TM\oplus L$ over $M^{2n}$ be endowed with 
the induced vector bundle connection denoted by $\hat\nabla$, 
that is,
$$
\hat\nabla_X(Y+\eta)=(\tilde\nabla_X(Y+\eta))_{TM\oplus L},
$$
where $\tilde\nabla$ is the Euclidean connection. It is easily
seen that the complex structure 
$\hat\J\in\Gamma(\text{Aut}(TM\oplus L))$ defined by 
$\hat\J(X+\eta)=JX+\J\eta$ is 
parallel, that is, we have
$$
(\nab_X\hat\J(Y+\eta))_{TM\oplus L}
=\hat\J((\nab_X(Y+\eta))_{TM\oplus L})
$$ 
for any $X\in\mathfrak{X}(M)$ and $\eta\in\Gamma(L)$.

\begin{proposition}\po\label{lambda} Let the vector subspaces  
$\mathcal{S}(x)=\Sal(\a|_{D\times D})(x)\subset L^{2\ell}(x)$  
satisfy that $\mathcal{S}(x)=L(x)$ for any $x\in M^{2n}$.
Then the $\hat\J$-invariant vector subbundle 
$\pi\colon\Lambda\to M^{2n}$ of $TM\oplus L^{2\ell}$ defined by
$$
\Lambda
=\spa\{(\nabla_ST)_{D^\perp}+\a(S,T)\colon S,T\in\Gamma(D)\}
$$ 
has rank $2\ell$. Moreover
$\Lambda\cap TM=0$ and $(\nabla_X\lambda)_{L^\perp}=0$ 
for $X\in\mathfrak{X}(M)$ and $\lambda\in\Gamma(\Lambda)$.
\end{proposition}

\proof This is Lemma $2.3$ of \cite{CD1}.\qed

\begin{theorem}\po\label{develop} Let 
$f\colon M^{2n}\to\R^{2n+p}$ be an embedded real Kaehler 
submanifold. Let $\pi\colon\hat\Lambda\to M^{2n}$ be a 
$\hat\J$-invariant vector subbundle of $TM\oplus L^{2\ell}$ 
of rank $2\ell$ satisfying the conditions 
$\hat\Lambda\cap TM=0$ and $(\nabla_X\lambda)_{L^\perp}=0$ 
for any $X\in\mathfrak{X}(M)$ and $\lambda\in\Gamma(\hat\Lambda)$. 
Let $N^{2n+2\ell}$ be an open neighborhood of the $0$-section 
$j\colon M^{2n}\to N^{2n+2\ell}$ of $\hat\Lambda$ such that 
the map $F\colon N^{2n+2\ell}\to\R^{2n+p}$ given by
\be\label{theextension}
F(\lambda)=f(\pi(\lambda))+\lambda
\ee
is an embedding. Then $F$ is a Kaehler extension of $f$ whose
second fundamental form satisfies
$\mathcal{N}_c(\a^F)=D\oplus\hat\Lambda$ at any point.
Moreover, the extension $F$ is minimal if and only if $f$ is
minimal.
\end{theorem}

\proof While the statement is slightly more general than that
of Theorem $2.5$ in \cite{CD1},  the proof remains unaltered, 
presented verbatim.\qed

\subsection{Complex ruled submanifold}

In the sequel, we assume that $L^{2\ell}$ splits as the 
orthogonal sum of vector subbundles 
$L^{2\ell}=\Sal\oplus{\cal R}$ where $\Sal=\Sal(\a|_{D\times D})$. 
From $\Sal^\perp={\cal R}\oplus L^\perp$ we have 
$\mathcal{N}_c(\a_{\Sal^\perp})
=D\cap\mathcal{N}_c(\a_{\mathcal{R}})$. 
Since  $\J\a(S,T)=\a(S,JT)\in\Sal$ for any $S,T\in\Gamma(D)$
then $\Sal$ and ${\cal R}$ are $\J$-invariant normal
subbundles and hence of even ranks. \vspace{1ex}

By Proposition $2.1$ in \cite{CD1} the distribution $D$ 
is integrable, and hence $M^{2n}$ carries a holomorphic 
foliation denoted by ${\cal F}$. 
Let $i\colon\Sigma\to M^{2n}$ be the inclusion of the 
leaf $\Sigma$ in ${\cal F}$ through $x\in M^{2n}$ and 
$g\colon\Sigma\to\R^{2n+p}$ the isometric immersion 
$g=f\circ i$. Its second fundamental form is
\be\label{secondfundg}
\a^g(y)(S,T)=f_*(i(y))(\nabla_{i_*S}i_*T)_{D^\perp}
+\a(i(y))(i_*S,i_*T)
\ee
for any $S,T\in\mathfrak{X}(\Sigma)$. Then the shape 
operators are
\be\label{shapeopg}
i_*A^g_{\xi\circ i}T=(A_\xi^fi_*T)_{D}\;\,\mbox{and}\;\,
i_*A^g_{(f_*X)\circ i}T=-(\nabla_{i_*T}X)_D
\ee 
for any $T\in\mathfrak{X}(\Sigma)$, $\xi\in\Gamma(L)$ 
and $X\in\Gamma(D^\perp)$. 

We have that $g(\Sigma)\subset f_*T_xM\oplus L(x)$. To see 
this, just observe that the normal bundle of $g$ splits 
orthogonally as 
$N_g\Sigma=i^*(f_*D^\perp\oplus L\oplus L^\perp)$, that
$i_*N_1^g\subset f_*   D^\perp\oplus L$ and that 
the condition $(\mathcal{C}_2$) gives that the vector subbundle 
$i^*L^\perp$ of $N_g\Sigma$ is constant in $\R^{2n+p}$. Therefore, 
when it is convenient we will consider $g(\Sigma)$ as a submanifold 
of $\R^{2n+2\ell}=f_*T_xM\oplus L^{2\ell}(x)\subset\R^{2n+p}$.

\begin{proposition}\po\label{charD2} For any $\Sigma\in {\cal F}$
the submanifold $g\colon\Sigma\to\R^{2n+2\ell}$ is holomorphic 
and the map $\psi\colon\Sal(\a^g)\to\Sal(\a|_{D\times D})$ 
defined by
\be\label{defder}
\psi(\a^g(S,T))=\a(i_*S,i_*T)
\ee
is an isomorphism. 
\end{proposition}

\proof This is Proposition $2.2$ in \cite{CD1}.\qed

\begin{theorem}\po\label{ell} Let $f\colon M^{2n}\to\R^{2n+p}$ 
be a real Kaehler submanifold that satisfies that 
$\dim N_1^f(x)=m\leq n-1$ 
and $\nu^c_f(x)<2n-2m$ at any $x\in M^{2n}$.  
If $f$ restricted to any open subset of 
$M^{2n}$ does not admit the Kaehler extension 
provided by \eqref{theextension}, then the following 
conclusions hold: 
\vspace{1ex}

\noindent $(I)$ If $\rank L=2$ then $d<2n-m$ and $f$ is 
$D$-complex ruled.
\vspace{1ex}

\noindent $(I\!I)$ If $\rank L=4$ then we have:
\begin{itemize}
\item[(a)] If $\Sal=0$ then $f$ is $D$-complex ruled with $3d<6n-2m$.
\item[(b)] If $\rank\Sal=2$, $d>2$ and none of the subbundles $N_1^g$ 
of the holomorphic submanifolds $g\colon\Sigma^d\to\R^{2n+2\ell}$ are 
parallel then $3d< 6n-2m+2$ and $f$ is $D'$-complex ruled where 
$D'\subset D$ satisfies $\rank D'\geq d-2$.
\item[(c)] If $\rank\Sal=2$, $d>n$ and the $N_1^g$'s 
are all parallel then $2d<4n-m$ and $f$ is $D'$-complex ruled 
where $D'\subset D$ satisfies $\rank D'\geq 2d-2n$.
\end{itemize}

\noindent $(I\!I\!I)$ If $\rank L=6$ then we have:
\begin{itemize}
\item[(a)] If $\Sal=0$ then $f$ is $D$-complex ruled with $2d<4n-m$ .
\item[(b)] If $\rank\Sal=2$, $d>2$ and none of the subbundles $N_1^g$ 
of the holomorphic submanifolds $g\colon\Sigma^d\to\R^{2n+2\ell}$ are 
parallel then $2d<4n-m+1$ and $f$ is $D'$-complex ruled where 
$D'\subset D$ satisfies $\rank D'\geq d-2$.
\item[(c)] If $\rank\Sal=2$, $3d>4n$ and the $N_1^g$'s  are all 
parallel then $3d<6n-m$ and $f$ is $D'$-complex ruled where 
$D'\subset D$ and $\rank D'\geq 3d-4n$.
\item[(d)] If $\rank\Sal=4$, $3d>4n+2$ and none of the 
$N_1^g$'s are parallel then $3d< 6n-m+2$ and $f$ is $D'$-complex 
ruled where $D'\subset D$ satisfies $\rank D'\geq 3d-4n-2$.
\item[(e)] If $\rank\Sal=4$, $2d>3n$ and the $N_1^g$'s are all 
parallel then $4d<8n-m$ and $f$ is $D'$-complex ruled where 
$D'\subset D$ satisfies $\rank D'\geq 4d-6n$.
\end{itemize}
Moreover, in all of the above cases in which $f$ is $D'$-complex 
ruled, it should be understood this to hold on connected components of 
an open dense subset of $M^{2n}$.
\end{theorem}

For the proof of the above result we require a series 
of lemmas given next.

\begin{lemma}\po\label{satisfies} Assume that for any 
$\Sigma\in{\cal F}$ the normal subbundle $N_1^g$ of 
$g\colon\Sigma\to\R^{2n+2\ell}$ is parallel in the 
normal connection. Then also $\Sal$ satisfies the 
conditions $(\mathcal{C}_1)$ and $(\mathcal{C}_2)$.
\end{lemma}

\proof Just using that $\Sal$ is $\J$-invariant
it is easily seen that it satisfies the condition 
$(\mathcal{C}_1)$. Hence we argue for the condition 
$(\mathcal{C}_2)$. We obtain from \eqref{shapeopg} that 
$i_*A^g_{\xi\circ i}T=(A^f_\xi i_*T)_D=0$ for any 
$\xi\in\Gamma({\cal R})$ and $T\in\mathfrak{X}(\Sigma)$. 
Thus $i^*{\cal R}\subset (N_1^g)^\perp$. Now using that 
$N_1^g{}^\perp$ is parallel by assumption and 
\eqref{secondfundg}, we obtain
\begin{align*}
0&=\<\a^g(S,T),{}^g\nabla_Z^\perp(\xi\circ i)\>
=\<f_*(\nabla_{i_*S}i_*T)_{D^\perp}
+\a(i_*S,i_*T),-f_*(A_{\xi\circ i}^fi_*Z)_{D^\perp}
+{}^f\nabla_{Z}^\perp\xi\>\\
&=\<\a(i_*S,i_*T),{}^f\nabla_{Z}^\perp\xi\>
\end{align*}
for any $S,T\in\mathfrak{X}(\Sigma)$, 
$Z\in i^*\mathcal{N}_c(\a_{\Sal^\perp})\subset\mathfrak{X}(\Sigma)$
and $\xi\in\Gamma({\cal R})$. 
This and the condition $(\mathcal{C}_2)$ for $L$ give 
that $\Sal$ is parallel along 
$\mathcal{N}_c(\a_{\Sal^\perp})$ in the normal 
connection of $f$.\qed 

\begin{lemma}\po\label{intkernel}
Let ${\cal L}^{2r}\subset L^{2\ell}$ be a $\J$-invariant 
vector subbundle. If the $J$-invariant vector subbundle 
$E^e\subset D$ satisfies 
$A_\xi E\subset D^\perp$ for any  $\xi\in\Gamma({\cal L})$,  
then $\nu^c(\a_{{\cal L}\oplus L^\perp})\geq e+rd-2nr$.
\end{lemma}

\proof If $\mathcal{{\cal L}}=\spa\{\eta_j,\J\eta_j,\,1\leq j\leq r\}$ 
and since $A_{\J\eta}=J\circ A_\eta$ for any $\eta\in\Gamma({\cal L})$, 
then
\be\label{intkerneleq}
\cap_{j=1}^r\ker A_{\eta_j}\cap E
\subset\cap_{j=1}^r\ker A_{\eta_j}\cap D
=\mathcal{N}_c(\a_{\cal L})\cap\mathcal{N}_c(\a_{L^\perp})
=\mathcal{N}_c(\a_{{\cal L}\oplus L^\perp}).
\ee 
Hence 
$\nu^c(\a_{{\cal L}\oplus L^\perp})
\geq\dim\cap_{j=1}^r\ker A_{\eta_j}\cap E$.
From $A_{\eta_1}|_{E}\colon E^e\to D^\perp$  we obtain 
$$
\dim\ker A_{\eta_1}\cap E\geq e+d-2n.
$$
If $r\geq 2$ from $A_{\eta_2}|_{\ker A_{\eta_1}\cap E}
\colon\ker A_{\eta_1}\cap E\to D^\perp$ we have
$\dim\cap_{j=1}^2\ker A_{\xi_j}\cap E\geq e+2d-4n$
and, similarly, that
$\dim\cap_{j=1}^r\ker A_{\xi_j}\cap E\geq e+rd-2nr$.\qed

\begin{lemma}\po\label{cases} The following facts hold:
\begin{itemize} 
\item[(i)] If ${\cal R}\neq 0$ and $0\neq\xi\in\Gamma({\cal R})$ 
then
\be\label{c:snullity}
\nu^c(\a_{\spa\{\xi,\J\xi\}\oplus L^\perp})
\geq\max\{0,2d-2n\}.
\ee 
\item[(ii)] If $\Sal=0$ then $f$ is $D$-complex ruled
and 
\be\label {c:snullity2}
\nu^c_f\geq\max\{0,(\ell+1)d-2n\ell\}.
\ee
\end{itemize}
\end{lemma}

\proof Since $\mathcal{R}$ is $\J$-invariant then 
$\mathcal{L}=\spa\{\xi,\J\xi\}$ is $\J$-invariant. 
Having that $A_\eta D\subset D^\perp$ 
for any $\eta\in\Gamma(\mathcal{L})$ then \eqref{c:snullity} 
follows from Lemma \ref{intkernel}.

From the Codazzi equation $(\nabla_X^\perp\a)(S,T)
=(\nabla_S^\perp\a)(X,T)$ and the condition $(\mathcal{C}_2)$ 
we obtain that
\be\label{codhojapluri}
\a_{L^\perp}(X,\nabla_ST)+(\nabla_X^\perp\a(S,T))_{L^\perp}=0
\ee 
for any $S,T\in\Gamma(D)$ and $X\in\mathfrak{X}(M)$.  
If $\Sal=0$, then the distribution $D$ is a totally 
geodesic in $M^{2n}$. Hence, its leaves are totally geodesic 
submanifolds in the ambient space. Since $A_\xi D\subset D^\perp$ 
if $\xi\in\Gamma(L)$ then Lemma \ref{intkernel} 
applied to $L^{2\ell}$ gives \eqref{c:snullity2}.\qed

\begin{lemma}\po\label{pnucleohojas}
Let $\E\subset D$ be the vector subbundle given by 
$$
\E=\spa\{T\in\Gamma(D)\colon (\nabla_ST)_{D^\perp}+\a(S,T)=0
\;\,\mbox{for any}\;\,S\in\Gamma(D)\}.
$$ 
If $\E\neq 0$, then it is $J$-invariant and $f$ is $\E$-complex ruled. 
\end{lemma}

\proof Since $D$ is $J$-invariant then 
$J(\nabla_ST)_{D^\perp}=(\nabla_SJT)_{D^\perp}$
for any $S,T\in\Gamma(D)$. Hence
$$
\hat\J((\nabla_ST)_{D^\perp}+\a(S,T))
=(\nabla_SJT)_{D^\perp}+\a(S,JT)
$$
for any $S,T\in\Gamma(D)$. Thus $\E$ is $J$-invariant. 

From \eqref{secondfundg} for $\Sigma\in{\cal F}$ the holomorphic 
submanifold $g=f\circ i\colon\Sigma\to\R^{2n+p}$ satisfies  
\be\label{sigma}
\E|_\Sigma=i_*\mathcal{N}_c(\a^g).
\ee
We have that $\nabla_ST\in\Gamma(D)$  for any $T\in\Gamma(\E)$ 
and $S\in\Gamma(D)$. By \eqref{sigma} if $S\in\Gamma(\E)$ there 
is $\bar{S}\in\mathcal{N}_c(\a^g)$ such that $i_*\bar S=S|_\Sigma$.
Since the distribution $\mathcal{N}_c(\a^g)$ is totally geodesic 
then $^\Sigma\nabla_{\bar S}\bar{T}\in\mathcal{N}_c(\a^g)$ for any 
$S,T\in\Gamma(\E)$. From 
$(\nabla_Si_*\bar{T})_D=i_*{}^\Sigma\nabla_{\bar{S}}\bar{T}$
for any $S,T\in\Gamma(D)$ it follows that $\nabla_ST\in\Gamma(\E)$ 
for any  $S,T\in\Gamma(\E)$. Thus $\E$ is a totally geodesic 
distribution.  Since $\a(S,T)=0$ for any $T\in\Gamma(\E)$ and 
$S\in\Gamma(D)$ then $f$ is  $\E$-complex ruled.
\vspace{2ex}\qed

\noindent{\em Proof of Theorem \ref{ell}:} We prove part $(I)$. 
Proposition \ref{lambda} and Theorem \ref{develop} yield 
that $\Sal=0$. Then \eqref{c:snullity2} gives that $f$ is 
$D$-complex ruled and that $\nu^c_f\geq 2d-2n$. By assumption 
$\nu^c_f<2n-2m$, and hence $d<2n-m$.
\vspace{1ex}

The arguments for the proof of part $(I\!I)$ are contained in 
the more difficult arguments for $(a),(b)$ and $(c)$ in part 
$(I\!I\!I)$, and thus we only argue for the latter part. Then 
$\rank L=6$ and it follows from Proposition \ref{lambda} and 
Theorem \ref{develop} that $\rank\Sal=0,2$ or $4$.
\vspace{1ex} 

We prove $(a)$. We have from  \eqref{c:snullity2} that $f$ is 
$D$-complex ruled and $\nu^c_f\geq 4d-6n$. Since $\nu^c_f<2n-2m$ 
then $2d<4n-m$.  
\vspace{1ex}

We prove $(b)$. Given $\Sigma\in{\cal F}$, by assumption there are 
$S_0\in\mathfrak{X}(\Sigma)$ and $\eta\in\Gamma(N_1^g{}^\perp)$ 
such that $(^g\nabla_{S_0}^\perp\eta)_{N_1^g}\neq 0$. 
Proposition \ref{charD2} yields $\rank N_1^g=2$. The Codazzi 
equation gives 
$$
A^g_{(^g\nabla_S^\perp\eta)_{N_1^g}}T
=A^g_{(^g\nabla_T^\perp\eta)_{N_1^g}}S
$$
and, since $g$ is holomorphic, that $\nu^c(\a^g)\geq d-2>0$. 
Hence \eqref{sigma} yields $\rank\E\geq d-2$ and then 
Lemma \ref{pnucleohojas} that $f$ is $D'$-complex ruled 
where $D'=\E$.

From \eqref{shapeopg} and \eqref{sigma} we have 
$(A_\xi^fS)_D=i_*A^g_{\xi\circ i}\bar{S}=0$
for any $S\in\Gamma(\E)$ and $\xi\in\Gamma(L)$. Thus
$A_\xi^f\E\subset D^\perp$ for any $\xi\in\Gamma(L)$.
Then Lemma~\ref{intkernel} applied to the pair $\{L,\E\}$
yields $\nu^c_f\geq\rank\E+3d-6n$. Since $\rank\E\geq d-2$
and $\nu^c_f<2n-2m$ then $2d<4n-m+1$.
\vspace{1ex}

We prove $(c)$. Lemma \ref{satisfies} gives that $\Sal$ 
satisfies the conditions $(\mathcal{C}_1)$ and $(\mathcal{C}_2)$. 
Then Lemma~\ref{intkernel} applied to $\{{\cal R},D\}$ 
gives $\nu^c(\a_{\Sal^\perp}^f)\geq 3d-4n>0$.
Then part $(I)$ yields that $f$ is $D'$-complex ruled
where $D'=\mathcal{N}_c(\a_{\Sal^\perp}^f)$ and $d'<2n-m$. 
Hence $3d< 6n-m$.
\vspace{1ex}

We prove $(d)$. Assume first that $\nu^c_g\geq d-6$ for any 
$\Sigma\in{\cal F}$. As seen above, we have $A^f_\xi\E\subset D^\perp$
and then Lemma \ref{intkernel} applied to $\{L,{\E}\}$ gives
$\nu^c_f\geq 4d-6n-6$. This and $\nu^c_f<2n-2m$
yield $2d< 4n-m+3$. Since $m\geq 2\ell=6$ we have $d\leq 2n-2$.
Then $3d<6n-m+1$ and $d-6\geq3d-4n-2>0$. Hence, we have by 
\eqref{sigma} that $\rank\E\geq d-6\geq 3d-4n-2>0$. Finally,
Lemma \ref{pnucleohojas} gives that $f$ is $D'$-complex ruled 
where $D'=\E$ and $d'\geq 3d-4n-2$.

We now assume $\nu^c_g\leq d-8$. By assumption the tensor
$\phi\colon\mathfrak{X}(\Sigma)\times\Gamma(N_1^g{}^\perp)
\to\Gamma(N_1^g)$ given by 
$\phi(T,\eta)=({}^g\nabla_T^\perp\eta)_{N_1^g}$ is nontrivial. 
The Codazzi equation yields 
\be\label{codazziell3}
A_{\phi(S,\eta)}^gT=A_{\phi(T,\eta)}^gS
\ee
for $\eta\in\Gamma(N_1^g{}^\perp)$ and $S,T\in\mathfrak{X}(\Sigma)$.
Since $g$ is holomorphic in $\R^{2n+2\ell}$ then $\J$ 
induces an almost complex structure $\J_N$ on $N_1^g$.
We claim that $N_1^g\neq\Sal(\phi)+\J_N\Sal(\phi)$, 
thus assume otherwise. First suppose that there exists 
$\eta\in\Gamma(N_1^g{}^\perp)$ such that 
$N_1^g=\Ima\phi_\eta+\J_N\Ima\phi_\eta$. By \eqref{codazziell3}
we have  $\ker\phi_\eta\subset\ker A^g_{Im\phi_\eta}$. Since  
$A_{\J_N\delta}^g=J|_{\Sigma}\circ A_\delta^g$ holds if 
$\delta\in N_1^g$, then $\ker\phi_\eta\subset\ker A^g_{Im \phi_\eta}
\cap\ker A^g_{Im \J_N\phi_\eta}$. Having that Proposition \ref{charD2} 
gives $\rank N_1^g=\rank\Sal=4$, we obtain from 
$$
\mathcal{N}(\a^g)=\cap_{\delta\in N_1^g}\ker A^g_\delta
=\ker A^g_{Im \phi_\eta} \cap\ker A^g_{Im \J_N\phi_\eta}
$$
that $\nu^c_g\geq\dim\ker\phi_\eta\geq d-4$, which is 
a contradiction.

The other possible situation is that there are vector fields 
$\eta_j\in\Gamma(N_1^g{}^\perp)$, $j=1,2$, such that
$N_1^g=\sum_{j=1}^2(\Ima\phi_{\eta_j}+\J_N\Ima\phi_{\eta_j})$
where $1\leq\dim\Ima\phi_{\eta_j}\leq 2$ for $j=1,2$.   
Then $\ker\phi_{\eta_j}\subset\ker A^g_{Im \phi_{\eta_j}}
\cap\ker A^g_{Im \J_N\phi_{\eta_j}}$, $j=1,2$. 
Since 
$$
\mathcal{N}(\a^g)=\cap_{\delta\in N_1^g}\ker A^g_\delta=
\cap_{j=1}^2 \ker A^g_{Im \phi_{\eta_j}}
\cap\ker A^g_{Im \J_N\phi_{\eta_j}}
$$
hence $\nu^c_g\geq\dim\ker\phi_{\eta_1}\cap\ker\phi_{\eta_2}\geq d-4$.
But this is a contradiction proving the claim.

From the claim, there is an orthogonal decomposition 
$N_1^g=U^2\oplus V^2$ where $V^2=\Sal(\phi)+\J_N\Sal(\phi)$.
Hence there exists $\eta\in\Gamma(N_1^g{}^\perp)$ such that 
$V^2=\Ima\phi_\eta+\J_N\Ima\phi_\eta$. 
By \eqref{codazziell3} and since $g$ is holomorphic then 
$\ker\phi_\eta\subset\ker A^g_{Im \phi_\eta}
\cap\ker A^g_{Im \J_N\phi_\eta}$.  It follows that
$$
{\cal D}_1=\mathcal{N}(\a^g_{V^2})=\cap_{\delta\in V^2}\ker A^g_\delta=
\ker A^g_{Im \phi_\eta} \cap\ker A^g_{Im \J_N\phi_\eta}
$$
and thus $\rank{\cal D}_1\geq\dim\ker\phi_\eta\geq d-2$.

If  ${\cal D}_1=D$ then $N_1^g=U^2$ contradicting that 
$\rank N_1^g=\rank\Sal=4$.
Hence  ${\cal D}_1\neq D$, and since $g$ is holomorphic then 
${\cal D}_1$ is $J|_\Sigma$-invariant and hence 
$\rank{\cal D}_1=d-2$.

We argue that $U^2$ is parallel along ${\cal D}_1$, thus  
assume otherwise. The Codazzi equation gives
$$
A^g_{(^g\nabla_S^\perp\xi)_{U^2}}T-
A^g_{(^g\nabla_T^\perp\xi)_{U^2}}S
=A_\xi^g\,[S,T]\in\Gamma({\cal D}_1^\perp)
$$
for any $S,T\in\Gamma({\cal D}_1)$ and $\xi\in\Gamma(U^\perp)$. 
Since $\rank{\cal D}_1^\perp=2$ it follows that $\nu^c_g\geq d-6$,
and this is a contradiction.

Since we assumed that $\nu^c_g\leq d-8$ then 
${\cal D}_1\neq\mathcal{N}(\a^g)$. 
It follows that
\be\label{ell3t2}
U^2=\spa\{\a^g(S,T)\colon S\in\Gamma(D)
\;\mbox{and}\;T\in\Gamma({\cal D}_1)\}.
\ee 

Let $L=\psi(U^2)\oplus W$ be an orthogonal decomposition where 
$\psi\colon N_1^g\to\Sal$ is the isomorphism given by \eqref{defder}. 
Since ${\cal D}_1$ is $J|_{\Sigma}$-invariant then $\psi(U^2)$ is 
$\J$-invariant. For any $\xi\in\Gamma(W)$ we have that
$\xi\circ i\in\Gamma(U^2{}^\perp)$ where
$U^2{}^\perp=V^2\oplus N_1^g{}^\perp$. From \eqref{shapeopg} we
obtain that $A_\xi^f i_*{\cal D}_1\subset D^\perp$ for  
$\xi\in\Gamma(W)$. Let $W=\spa\{\xi_1,\J\xi_1,\xi_2,\J\xi_2\}$.
By \eqref{intkerneleq}, the estimate in Lemma \ref{intkernel} 
and using that $\rank {\cal D}_1=d-2$, there is 
$D'=\cap_{j=1}^2\ker A^f_{\xi_j}\cap i_*{\cal D}_1$
with $D'\subset\mathcal{N}_c(\a_{\psi(U)^\perp})$ and 
$d'\geq 3d-4n-2$. We prove that 
$D'=\mathcal{N}_c(\a_{\psi(U)^\perp})$. Assume otherwise. 
We know that $A_\xi^f i_*{\cal D}_1\subset D^\perp$  and  
$A_\xi^f\mathcal{N}_c(\a_{\psi(U)^\perp})=0$ for any 
$\xi\in\Gamma(W)$. Since $D'$ and 
$\mathcal{N}_c(\a_{\psi(U)^\perp})$ are both $J$-invariant
and $\dim{\cal D}_1=d-2$ holds, then
$D=i_*{\cal D}_1+\mathcal{N}_c(\a_{\psi(U)^\perp})$. Hence 
we have that $A_\xi^f D\subset D^\perp$ for any $\xi\in\Gamma(W)$.  
Thus $\rank\Sal(\a|_{D\times D})=\rank\Sal=2$, 
and this is a contradiction.  

Since $U^2$ is parallel along ${\cal D}_1$ and
$D'=\mathcal{N}_c(\a_{\psi(U^2)^\perp})$ we obtain using
\eqref{secondfundg} that
\begin{align*}
0&=\<\a^g(S,T),{}^g\nabla_Z^\perp(\xi\circ i)\>
=\<f_*(\nabla_{i_*S}i_*T)_{D^\perp}
+\a(i_*S,i_*T),-f_*(A_{\xi\circ i}^fi_*Z)_{D^\perp}
+{}^f\nabla_{i_*Z}^\perp\xi\>\\
&=\<\a(i_*S,i_*T),{}^f\nabla_{i_*Z}^\perp\xi\>
\end{align*}
for any $S\in\mathfrak{X}({\cal D}_1)$, $T\in\mathfrak{X}(\Sigma)$,
$Z\in\Gamma(i^*D' )$ and $\xi\in\Gamma(W)$.
By \eqref{ell3t2} this gives that $\psi(U^2)$ is parallel
along $D'$. Since $\psi(U^2)$ is $\J$-invariant then it
satisfies the conditions $(\mathcal{C}_1)$ and $(\mathcal{C}_2)$.
Hence part $(I)$ gives that $f$ is $D'$-complex ruled
with $d'<2n-m$. Since $d'\geq 3d-4n-2$ then $3d<6n-m+2$.
\vspace{1ex}

We prove $(e)$. Lemma \ref{intkernel} applied to 
$\{{\cal R},D\}$ gives that $D_0=\mathcal{N}_c(\a_{\Sal^\perp})$ 
has rank $d_0\geq 2d-2n$. By assumption $2d>3n$ and hence $d_0>n$.
Lemma \ref{satisfies} yields that $\Sal$ satisfies the conditions 
$(\mathcal{C}_1)$ and $(\mathcal{C}_2)$. Then part $(I\!I)$ applied 
to $\Sal$ as $L$ gives that $f$ is $D'$-complex ruled with 
$D'\subset D_0$ and that one of following possibilities holds: 
$d'=d_0\geq 2d-2n$,  $d'\geq d_0-2\geq 2d-2n-2$ or
$d'\geq 2d_0-2n\geq 4d-6n$. Since $m\geq 2\ell=6$ then part $(I\!I)$  
also gives that $2d_0<4n-m$. Hence $4d-4n\leq 2d_0< 4n-m$ and 
thus $4d<8n-m$ which gives $d\leq2n-2$. Hence 
$d'\geq\min\{2d-2n-2,4d-6n\}=4d-6n$.\qed

\section{A certain vector subbundle}
\label{s:complexvb}

In this section, associated to a real Kaehler submanifold 
we construct a normal vector subbundle that meets
the conditions $(\mathcal{C}_1)$ and $(\mathcal{C}_2)$ 
introduced in the preceding section.
\vspace{2ex}

Throughout this section and the following ones, we will 
define numerous vector subspaces on a pointwise basis, either 
as the images or as the kernels of specific tensor fields 
situated on the submanifold. Given that our goal is confined 
to local results, in order to prevent excessive repetition, 
it is understood, without the need for additional mention, 
that we are  working restricted to connected components within 
an open dense subset of the manifold where these vector subspaces 
maintain constant dimensions, thus forming smooth vector subbundles.
\vspace{1ex}

We first argue that Theorem \ref{alglemma} applies pointwise 
to the second fundamental form of a real Kaehler submanifold.
The following result will be of use in the remaining of the 
paper.

\begin{proposition}\po 
Let $f\colon M^{2n}\to\R^{2n+p}$
be a real Kaehler submanifold. Let the bilinear forms
$\gamma,\beta\colon T_xM\times T_xM\to N_1^f(x)\oplus N_1^f(x)$ 
be defined  by \eqref{gamma} and \eqref{beta} in terms of the 
second fundamental form $\a$ of $f$ at $x\in M^{2n}$. Then
both bilinear forms are flat and satisfy the condition 
\eqref{productflat2}.
\end{proposition}

\proof Since  the Riemannian curvature tensor of a Kaehler 
manifold satisfies that $R(X,Y)JZ=JR(X,Y)Z$ for any 
$X,Y,Z\in T_xM$, then a roughly short straightforward 
computation using the Gauss equation of $f$ gives the 
result.\vspace{1ex}\qed 

In the sequel, associated to a real Kaehler submanifold 
$f\colon M^{2n}\to\R^{2n+p}$ we construct vector subbundles 
$N_1^f\supset L_0\supset L_1\supset L_2=L$ 
and their associated  vector subbundles 
$TM\supset D_0\supset D_1\supset D_2=D$ defined by 
$D_j=\mathcal{N}_c(\a_{L_j^\perp})$, $j=0,1,2$. In addition,
it is shown that there are associated isometric complex 
structures $\J_j\in\Gamma(\text{Aut}(L_j))$, $0\leq j\leq 2$, 
which are parallel in the induced connections on $L_j$ for 
$j=1,2$ and satisfy that 
$\J_j\a_{L_j}(X,Y)=\a_{L_j}(X,JY),0\leq j\leq 2$,
for any $X,Y\in\mathfrak{X}(M)$. 
\vspace{2ex}

Assume that $\Sal(\gamma)$ is a degenerate vector 
subspace. Part $(iv)$ of Proposition~\ref{even} yields that
$\Omega=\Sal(\a_\Omega)$, where 
$0\neq\Omega=\pi_1(\Sal(\gamma)\cap\Sal(\gamma)^\perp)$ 
and $\pi_1$ denotes taking the first component.
Let $N_1^f(x)=\Omega\oplus P$ be an orthogonal decomposition.  
Then the bilinear form 
$\gamma_P\colon T_xM\times T_xM\to P\oplus P$ defined by 
$\gamma_P=\pi_{P\times P}\circ\gamma$ is flat and 
$\Sal(\gamma_P)$ is a nondegenerate vector subspace. 

Let $D_0^{d_0}=\mathcal{N}_c(\a_P)$ and let $L_0\subset\Omega$ 
be the vector subspace  $L_0=\Sal(\a|_{T_xM\times D_0})$. 
Then
\be\label{condkerd0}
D_0^{d_0}=\mathcal{N}_c(\a_{L_0^\perp})
=\cap_{\xi\in L_0 ^\perp}\{\kerl A_\xi\cap\kerl A_\xi\circ J\}.
\ee
Let $\mathcal{L}\in\text{Aut}(\Omega)$ be given by
\be\label{condd0}
\mathcal{L}\a_\Omega(X,Y)=\a_\Omega(X,JY)
\;\,\mbox{for any}\;\,X,Y\in T_xM.
\ee
That $D_0$ is $J$-invariant yields that $L_0$ 
is $\mathcal{L}$-invariant, and thus $\dim L_0=2\ell_0$. 
\vspace{1ex}

Part $(iv)$ of Proposition \ref{even} gives that $\gamma_\Omega$
is a null bilinear form. Thus $\mathcal{L}\in\text{Aut}(\Omega)$ 
is an isometric complex structure. Since $L_0$ is 
$\mathcal{L}$-invariant, the vector bundle isometry 
$\J_0=\mathcal{L}|_{L_0}\in\Gamma(\text{Aut}(L_0))$ is an 
isometric complex structure that satisfies 
\be\label{condd00}
\J_0\a_{L_0}(X,Y)=\a_{L_0}(X,JY)\;\,\mbox{for any}\;\,X,Y\in T_xM.
\ee
Equivalently, we have 
$A_{\J_0\delta}=J\circ A_\delta=-A_\delta\circ J$ for any 
$\delta\in\Gamma(L_0)$.

\begin{lemma}\po\label{charK} The skew-symmetric tensor 
$\mathcal{K}(X)\in\Gamma(\text{End}(L_0))$  defined by
\be\label{kappa}
\mathcal{K}(X)\eta
=\J_0(\nabla_X^\perp\eta)_{L_0}-(\nabla_X^\perp\J_0\eta)_{L_0}
\ee
for any $X\in\mathfrak{X}(M)$ satisfies:
\begin{enumerate}
\item[(i)] $\mathcal{K}(X)\circ\J_0=-\J_0\circ\mathcal{K}(X)$.
\item[(ii)] $\mathcal{K}(Z)=0$ if $Z\in\Gamma(D_0)$.
\item[(iii)] $\mathcal{K}(X)\a(Y,Z)=\mathcal{K}(Y)\a(X,Z)$ 
if $X,Y\in\Gamma(D_0)$ or $Z\in\Gamma(D_0)$.
\item[(iv)] $\<\mathcal{K}(X)\a(Y,Z),\a(T,Z)\>=0$ if 
$X,Y,T\in\mathfrak{X}(M)$
and $Z\in\Gamma(D_0)$.
\end{enumerate}
\end{lemma}

\proof Part $(i)$ follows using that $\J_0^2=-I$. 
Applying $\J_0$ to the $L_0$-component of the Codazzi equation 
$(\nabla_X^\perp\a)(Y,Z)=(\nabla_Y^\perp\a)(X,Z)$, subtracting 
the $L_0$-component of the equation
$(\nabla_X^\perp\a)(Y,JZ)=(\nabla_Y^\perp\a)(X,JZ)$ and then 
using \eqref{condd0} give
$$
\J_0(\nabla_X^\perp\a(Y,Z))_{L_0}-(\nabla_X^\perp\a(Y,JZ))_{L_0}
=\J_0(\nabla_Y^\perp\a(X,Z))_{L_0}-(\nabla_Y^\perp\a(X,JZ))_{L_0}
$$
for any $X,Y,Z\in\mathfrak{X}(M)$. Then part $(iii)$ follows 
using \eqref{condkerd0} and \eqref{condd00}.

We denote
$$
(X_1,X_2,X_3,X_4,X_5)=
\<\mathcal{K}(X_1)\a_{L_0}(X_2,X_3),\a_{L_0}(X_4,X_5)\>
$$
for any $X_1,X_2,X_3,X_4,X_5\in\mathfrak{X}(M)$.  The
skew-symmetry of $\mathcal{K}(X)$ and part $(iii)$ yield
\begin{align*}
(Z_1,Y,Z_2,X,Z_3)&=(Y,Z_1,Z_2,Z_3,X)
=-(Y,Z_3,X,Z_1,Z_2)=-(X,Z_3,Y,Z_1,Z_2)\\
&=(X,Z_1,Z_2,Z_3,Y)=(Z_2,Z_1,X,Z_3,Y)=-(Z_2,Z_3,Y,Z_1,X)\\
&=-(Z_3,Z_2,Y,Z_1,X)=(Z_3,Z_1,X,Z_2,Y)=(Z_1,Z_3,X,Z_2,Y)\\
&=-(Z_1,Z_2,Y,Z_3,X)=-(Z_1,Y,Z_2,X,Z_3)=0
\end{align*}
for any $Z_1,Z_2,Z_3\in\Gamma(D_0)$ and $X,Y\in\mathfrak{X}(M)$.  
This gives part $(ii)$.

Finally, part $(iv)$ follows from
\begin{align*}
(X,Y,Z,T,Z)&=(Y,X,Z,T,Z)=-(Y,T,Z,X,Z)=-(T,Y,Z,X,Z)\\
&=(T,X,Z,Y,Z)=(X,T,Z,Y,Z)=-(X,Y,Z,T,Z)=0	
\end{align*}
for any $X,Y,T\in\mathfrak{X}(M)$ and $Z\in\Gamma(D_0)$.
\qed 
\vspace{1ex}

Let $L_1\subset L_0$ be the vector subbundle 
$L_1=\cap_{X\in\mathfrak{X}(M)}\ker\mathcal{K}(X)$ and let
$D_1\subset D_0$ be given by $D_1^{d_1}=\mathcal{N}_c(\a_{L_1^\perp})$. 
Then \eqref{kappa} gives
\be\label{parallell0}
(\nap_X\J_0\delta)_{L_0}=\J_0(\nap_X\delta)_{L_0}
\ee  
for any $X\in\mathfrak{X}(M)$ and $\delta\in\Gamma(L_1)$. 
Since $L_1$ is $\J_0$-invariant by part $(i)$ of 
Lemma \ref{charK}, then the vector bundle isometry 
$\J_1=\J_0|_{L_1}\in\Gamma(\text{Aut}(L_1))$ 
is a complex structure, and thus $\dim L_1=2\ell_1$.

If $L_0=L_1\oplus R_1$ is an orthogonal decomposition, 
we have by \eqref{condd00} that 
$$
\J_1\a_{L_1}(X,Y)+\J_0\a_{R_1}(X,Y)=\J_0\a_{L_0}(X,Y)
=\a_{L_0}(X,JY)=\a_{L_1}(X,JY)+\a_{R_1}(X,JY)
$$
for any $X,Y\in\mathfrak{X}(M)$. Therefore 
\be\label{condd1}
\J_1\a_{L_1}(X,Y)=\a_{L_1}(X,JY)\;\;\mbox{and}\;\;
\J_0\a_{R_1}(X,Y)=\a_{R_1}(X,JY)
\ee 
for any $X,Y\in\mathfrak{X}(M)$.  From \eqref{parallell0} we 
obtain that $\J_1$ is parallel, that is, that
\be\label{parallell1}
(\nap_X\J_1\eta)_{L_1}=\J_1(\nap_X\eta)_{L_1}
\ee 
for any $X\in\mathfrak{X}(M)$ and $\eta\in\Gamma(L_1)$.

\begin{lemma}\po\label{liebracketprop}
We have that $[Y,Z]\in\Gamma(D_0)$ if
$Y\in\Gamma(D_1)$ and $Z\in\Gamma(D_0)$.
\end{lemma}

\proof  If $\delta\in\Gamma(L_0)$, then the Codazzi equation 
$(\nabla_Z A)(\xi;Y)=(\nabla_Y A)(\xi;Z)$ for $\xi=\J_0\delta$ 
and that $A_{\J_0\delta}=J\circ A_\delta$ yield
$$
J(\nabla_Z A_\delta Y-A_\delta\nabla_ZY)
-A_{\nabla_Z^\perp\J_0\delta}Y
=J(\nabla_Y A_\delta Z-A_\delta\nabla_YZ)
-A_{\nabla_Y^\perp\J_0\delta}Z
$$
for any $Y,Z\in\mathfrak{X}(M)$.
Making use of the Codazzi equation for $\xi=\delta$, 
we obtain
$$
JA_{\nabla_Z^\perp\delta}Y-A_{\nabla_Z^\perp\J_0\delta}Y=
JA_{\nabla_Y^\perp\delta}Z-A_{\nabla_Y^\perp\J_0\delta}Z
\;\;\mbox{for any}\;\;Y,Z\in\mathfrak{X}(M).
$$
Hence 
$$
\<\a(Y,X),\nabla_Z^\perp\delta\>-\<\a(Y,JX),\nabla_Z^\perp\J_0\delta\>
=\<\a(Z,X),\nabla_Y^\perp\delta\>-\<\a(Z,JX),\nabla_Y^\perp\J_0\delta\>
$$
for any $\delta\in\Gamma(L_0)$ and $X,Y,Z\in\mathfrak{X}(M)$.  
Then \eqref{condd00}, \eqref{parallell0} and part $(ii)$ of 
Lemma \ref{charK} give
\be\label{ldtemp}
\<\a_{L_0^\perp}(Y,X),\nabla_{Z_1}^\perp\delta\>
-\<\a_{L_0^\perp}(Y,JX),\nabla_{Z_1}^\perp\J_0\delta\>=
\<\a(Z_1,JX),\mathcal{K}(Y)\delta\>
\ee
for any $Y,X\in\mathfrak{X}(M)$, $Z_1\in\Gamma(D_0)$ and 
$\delta\in\Gamma(L_0)$.  
Since $\delta=\a(Z,Z_2)\in\Gamma(L_1)$ if $Z_2\in\Gamma(D_1)$, 
we have using \eqref{condd00} that 

\be\label{ldtemp2}
\<\a_{L_0^\perp}(Y,X),\nabla_{Z_1}^\perp\a(Z,Z_2)\>-
\<\a_{L_0^\perp}(Y,JX),\nabla_{Z_1}^\perp\a(Z,JZ_2)\>=0
\ee
for any $Z_1\in\Gamma(D_0)$, $Z_2\in\Gamma(D_1)$ and 
$X,Y,Z\in\mathfrak{X}(M)$. 

On one hand, the Codazzi equation
$(\nap_{Z_1}\a)(Z,Z_2)=(\nap_{Z_2}\a)(Z,Z_1)$ 
and \eqref{condkerd0} give
\be\label{liebracket}
\a_{L_0^\perp}(Z,[Z_1,Z_2])
=(\nap_{Z_1}\a(Z,Z_2)-\nap_{Z_2}\a(Z,Z_1))_{L_0^\perp}
\ee
for any $Z_1,Z_2\in\Gamma(D_0)$ and $Z\in\mathfrak{X}(M)$.
On the other hand, the Codazzi equation 
$(\nap_{Z}\a)(Z_i,JZ_j)=(\nap_{Z_i}\a)(Z,JZ_j)$
and \eqref{condkerd0} yield
$$
(\nap_Z\a(Z_i,JZ_j))_{L_0^\perp}=
(\nap_{Z_i}\a(Z,JZ_j))_{L_0^\perp}-\a_{L_0^\perp}(Z,\nabla_{Z_i}JZ_j)
$$
for any $Z_i,Z_j\in\Gamma(D_0)$ and $Z\in\mathfrak{X}(M)$.
Since $\a(Z_1,JZ_2)=\a(JZ_1,Z_2)$ by \eqref{condd00}, then
\be\label{liebracket2}
\a_{L_0^\perp}(Z,J[Z_1,Z_2])=(\nap_{Z_1}\a(Z,JZ_2)
-\nap_{Z_2}\a(Z,JZ_1))_{L_0^\perp} 
\ee 
for any $Z_1,Z_2\in\Gamma(D_0)$ and $Z\in\mathfrak{X}(M)$.

First using \eqref{liebracket} and \eqref{liebracket2} 
and then \eqref{ldtemp2}, we obtain 
\begin{align}\label{new}
\<\a_{L_0^\perp}(Y,X)&,\a_{L_0^\perp}(Z,[Z_1,Z_2])\>-
\<\a_{L_0^\perp}(Y,JX),\a_{L_0^\perp}(Z,J[Z_1,Z_2])\>\nonumber\\
=&\,\<\a_{L_0^\perp}(Y,X),\nap_{Z_1}\a(Z,Z_2)\>
-\<\a_{L_0^\perp}(Y,X),\nap_{Z_2}\a(Z,Z_1)\>\nonumber\\
&-\<\a_{L_0^\perp}(Y,JX),\nap_{Z_1}\a(Z,JZ_2)\>
+\<\a_{L_0^\perp}(Y,JX),\nap_{Z_2}\a(Z,JZ_1)\>\nonumber\\
=&\,\<\a_{L_0^\perp}(Y,JX),\nap_{Z_2}\a(Z,JZ_1)\>
-\<\a_{L_0^\perp}(Y,X),\nap_{Z_2}\a(Z,Z_1)\>
\end{align}
for any $Z_1\in\Gamma(D_0)$, $Z_2\in\Gamma(D_1)$ and 
$X,Y,Z\in\mathfrak{X}(M)$. 
On the other hand, \eqref{ldtemp} gives 
$$
\<\a_{L_0^\perp}(Y,X),\nabla_{Z_2}^\perp\delta\>
-\<\a_{L_0^\perp}(Y,JX),\nabla_{Z_2}^\perp\J_0\delta\>=
\<\a(Z_2,JX),\mathcal{K}(Y)\delta\>
$$
for any $Z_2\in\Gamma(D_1)$, $\delta\in\Gamma(L_0)$ and 
$X,Y\in\mathfrak{X}(M)$. If $Z_1\in\Gamma(D_0)$ then we 
obtain for $\delta=\a(Z,Z_1)\in\Gamma(L_0)$ using 
\eqref{condd00} that 
\begin{align*}
\<\a_{L_0^\perp}(Y,X),&\nabla_{Z_2}^\perp\a(Z,Z_1)\>
-\<\a_{L_0^\perp}(Y,JX),\nabla_{Z_2}^\perp\a(Z,JZ_1)\>\\
&=\<\a(Z_2,JX),\mathcal{K}(Y)\a(Z,Z_1)\>
=-\<\mathcal{K}(Y)\a(Z_2,JX),\a(Z,Z_1)\>=0
\end{align*}
for any $Z_1\in\Gamma(D_0)$, $Z_2\in\Gamma(D_1)$ and 
$X,Y,Z\in\mathfrak{X}(M)$. Then \eqref{new} gives
\be\label{car}
\<\a_{L_0^\perp}(Y,X),\a_{L_0^\perp}(Z,[Z_1,Z_2])\>-
\<\a_{L_0^\perp}(Y,JX),\a_{L_0^\perp}(Z,J[Z_1,Z_2])\>=0.
\ee    

Since $L_0\subset\Omega$, we have from \eqref{condd0} that
\begin{align*}
\J_0\a_{L_0}(X,Y)
&+\mathcal{L}\a_{L_0^\perp\cap\Omega}(X,Y)
=\mathcal{L}\a_\Omega(X,Y)
=\a_\Omega(X,JY)\\
&\;=\a_{L_0}(X,JY)+\a_{L_0^\perp\cap\Omega}(X,JY)
\end{align*}
for any $X,Y\in\mathfrak{X}(M)$. We obtain from \eqref{condd00} that 
\be\label{above}
\mathcal{L}\a_{L_0^\perp\cap\Omega}(X,Y)
=\a_{L_0^\perp\cap\Omega}(X,JY)
\;\,\mbox{for any}\;\,X,Y\in\mathfrak{X}(M).
\ee
It follows from \eqref{car} and \eqref{above} that
$$
\lp\gamma_P(Y,X),\gamma_P(Z,[Z_1,Z_2])\rp=0
$$
for any $Z_1\in\Gamma(D_0)$, $Z_2\in\Gamma(D_1)$ and
$X,Y,Z\in\mathfrak{X}(M)$.  Being $\Sal(\gamma_P)$
nondegenerate, then $\gamma_P(Z,[Z_1,Z_2])=0$ for any
$Z\in\mathfrak{X}(M)$. From \eqref{car} for $JX$ 
instead of $X$, we have 
$$
\<\a_{L_0^\perp}(Y,X),\a_{L_0^\perp}(Z,J[Z_1,Z_2])\>+
\<\a_{L_0^\perp}(Y,JX),\a_{L_0^\perp}(Z,[Z_1,Z_2])\>=0
$$
for any $Z_1\in\Gamma(D_0)$, $Z_2\in\Gamma(D_1)$ and 
$X,Y,Z\in\mathfrak{X}(M)$. Then $\gamma_P(Z,J[Z_1,Z_2])=0$ 
for any $Z\in\mathfrak{X}(M)$, and therefore $[Z_1,Z_2]\in D_0$.
\vspace{2ex}\qed 

Let $L\subset L_1$ be the vector subbundle defined by
\be\label{defl2}
L=\{\delta\in\Gamma(L_1)\colon 
(\nabla_X^\perp\delta)_{L_0^\perp}=0
=(\nabla_X^\perp\J_1\delta)_{L_0^\perp} 
\;\mbox{for all}\;X\in\Gamma(D_0)\}
\ee 
and let $D\subset D_1$ be given by $D^d=\mathcal{N}_c(\a_{L^\perp})$.
Since $L$ is $\J_1$-invariant, then the vector bundle isometry 
$\J=\J_1|_L\in\Gamma(\text{Aut}(L))$ is a complex structure.   

\begin{proposition}\po The vector subbundle $L$ 
satisfies the conditions $(\mathcal{C}_1)$ and $(\mathcal{C}_2)$.
\end{proposition}

\proof We have from \eqref{condd1} and \eqref{parallell1} that 
$$
\J\a_L(X,Y)=\a_L(X,JY)\;\;\mbox{and}\;\;
\J(\nabla_X^\perp\delta)_L=(\nabla_X^\perp\J\delta)_L
$$
for any $X,Y\in\mathfrak{X}(M)$ and $\delta\in\Gamma(L)$.
Thus the condition $(\mathcal{C}_1)$ holds. 
\vspace{1ex}

The main ingredient to prove $(\mathcal{C}_2)$ 
is that $L_1$ is parallel along $D$, namely, that
\be\label{claim}
\nabla_Z^\perp\delta\in\Gamma(L_1)\;\,\mbox{for any}\;\,
Z\in\Gamma(D)\;\;\mbox{and}\;\;\delta\in\Gamma(L_1).
\ee 
If $L_0=L\oplus R$ is an orthogonal decomposition, 
that $L_0=\Sal(\a|_{TM\times D_0})$ gives
\be\label{ldimage}
R=\Sal(\a_R|_{TM\times D_0}).
\ee 
Lemma \ref{liebracketprop} and the Codazzi equations
$(\nabla_Y^\perp\a)(Z,X)=(\nabla_Z^\perp\a)(Y,X)$ yield 
$$
(\nabla_Y^\perp\a_{L_1}(Z,X))_{L_0^\perp}
=(\nabla_Z^\perp\a_{L_0}(Y,X))_{L_0^\perp}
$$
for any $Y\in\Gamma(D_0)$, $Z\in\Gamma(D_1)$ 
and $X\in\mathfrak{X}(M)$. Then \eqref{defl2} gives
$$
(\nabla_Y^\perp\a_{L^\perp\cap L_1}(Z,X))_{L_0^\perp}
=(\nabla_Z^\perp\a_R(Y,X))_{L_0^\perp}
$$
for any $Y\in\Gamma(D_0)$, $Z\in\Gamma(D_1)$ and 
$X\in\mathfrak{X}(M)$. Since the left-hand-side vanishes 
if $Z\in\Gamma(D)$, it follows from \eqref{defl2} and 
\eqref{ldimage} that $L_0$ is parallel along $D$, that is, 
that
\be\label{ldpara1}
\nabla_Z^\perp\delta\in \Gamma(L_0)\;\mbox{for all}\;
Z\in\Gamma(D)\;\;\mbox{and}\;\;\delta\in\Gamma(L_0).
\ee  

On one hand, if $L_0=L_1\oplus R_1$ is an orthogonal 
decomposition and since $\mathcal{K}(X)$ is skew-symmetric, 
then
\be\label{ldspanR}
R_1=\spa\{\mathcal{K}(X)\mu\colon \mbox{for any}\;X\in\mathfrak{X}(M)
\;\;\mbox{and}\;\;\mu\in\Gamma(R_1)\}.
\ee
On the other hand, by \eqref{condd00} we have that
$[A_{\J_0\delta},A_{\J_0\mu}]=[A_\delta,A_\mu]$ for any 
$\delta,\mu\in\Gamma(L_0)$. Then the Ricci equation gives
$$
\<R^\perp(X,Y)\delta,\mu\>=\<R^\perp(X,Y)\J_0\delta,\J_0\mu\>
$$
for any $X,Y\in\mathfrak{X}(M)$ and $\delta,\mu\in\Gamma(L_0)$.  
Now using \eqref{parallell0} we easily obtain
$$
\<\nabla_Y^\perp\delta,\nabla_X^\perp\mu\>
-\<\nabla_Y^\perp\J_0\delta,\nabla_X^\perp\J_0\mu\>=
\<\nabla_X^\perp\delta,\nabla_Y^\perp\mu\>-
\<\nabla_X^\perp\J_0\delta,\nabla_Y^\perp\J_0\mu\>
$$
for any $\delta\in\Gamma(L_1)$,  $\mu\in\Gamma(R_1)$
and $X,Y\in\mathfrak{X}(M)$. Hence, using \eqref{ldpara1} 
it follows that
\begin{align*}
\<\J_0\nabla_Y^\perp\delta,\J_0(\nabla_X^\perp\mu)_{L_0}\>
&-\<\nabla_Y^\perp\J_0\delta,(\nabla_X^\perp\J_0\mu)_{L_0}\>\\
&=\<\J_0(\nabla_X^\perp\delta)_{L_0},\J_0\nabla_Y^\perp\mu\>-
\<(\nabla_X^\perp\J_0\delta)_{L_0},\nabla_Y^\perp\J_0\mu\>
\end{align*}
for any $Y\in\Gamma(D)$, $\delta\in\Gamma(L_1)$,  
$\mu\in\Gamma(R_1)$ and $X\in\mathfrak{X}(M)$. Thus 
$$
\<\J_0\nabla_Y^\perp\delta,\mathcal{K}(X)\mu\>
=\<\J_0\nabla_X^\perp\delta,\mathcal{K}(Y)\mu\>=0
$$
for any $Y\in\Gamma(D)$, $\delta\in\Gamma(L_1)$, 
$\mu\in\Gamma(R_1)$ and $X\in\mathfrak{X}(M)$, where the first 
equality follows from \eqref{parallell0} and that $\delta\in\Gamma(L_1)$, 
and the second from part $(ii)$ of Lemma \ref{charK}. 
Then part $(i)$ of Lemma \ref{charK} and \eqref{ldspanR} yield
$\<\nabla_Y^\perp\delta,\J_0\mu\>=0$
for any $Y\in\Gamma(D)$, $\delta\in\Gamma(L_1)$ and
$\mu\in\Gamma(R_1)$. Since $R_1$ is 
$\J_0$-invariant by \eqref{condd1} then \eqref{claim} follows. 

We have 
$$
\<R^\perp(Y,Z)\delta,\xi\>=\<[A_\xi,A_\delta]Y,Z\>=0
$$
for any $Y,Z\in\Gamma(D_0)$, $\delta\in\Gamma(L_0)$ 
and $\xi\in\Gamma(L_0^\perp)$. Then \eqref{defl2}, 
\eqref{ldpara1} and Lemma \ref{liebracketprop} give
$$
0=\<R^\perp(Y,Z)\delta,\xi\>
=\<\nabla_Y^\perp\nabla_Z^\perp\delta,\xi\>
-\<\nabla_Z^\perp\nabla_Y^\perp\delta,\xi\>
-\<\nabla_{[Y,Z]}^\perp\delta,\xi\>
=\<\nabla_Y^\perp\nabla_Z^\perp\delta,\xi\>
$$
for any $Y\in\Gamma(D_0)$, $Z\in\Gamma(D)$, $\delta\in\Gamma(L)$ 
and $\xi\in\Gamma(L_0^\perp)$. Hence
$\nabla_Y^\perp\nabla_Z^\perp\delta\in\Gamma(L_0)$.
Since $L$ is $\J_1$-invariant we also have 
$\nabla_Y^\perp\nabla_Z^\perp\J_1\delta\in\Gamma(L_0)$. 
From \eqref{parallell1} and \eqref{claim} we obtain
$\nabla_Z^\perp\J_1\delta=\J_1(\nabla_Z^\perp\delta)$ and 
thus $\nabla_Y^\perp\J_1\nabla_Z^\perp\delta\in\Gamma(L_0)$.
Condition $(\mathcal{C}_2)$ is that 
$\nabla_Z^\perp\delta\in\Gamma(L)$ for any $Z\in\Gamma(D)$ and 
$\delta\in\Gamma(L)$. Since $\nabla_Z^\perp\delta\in\Gamma(L_1)$ 
by \eqref{claim}, then by \eqref{defl2} we need to have that 
$\nabla_Y^\perp\nabla_Z^\perp\delta$, 
$\nabla_Y^\perp\J_1\nabla_Z^\perp\delta\in\Gamma(L_0)$, 
which was already proved.\qed

\section{Some estimates}

This section is dedicated to furnishing estimates for 
the ranks of diverse vector subbundles that have been 
introduced in the previous section.
\vspace{2ex}

We assume in the sequel that the real Kaehler submanifold 
$f\colon M^{2n}\to\R^{2n+p}$ with codimension $p\leq n-1$ 
satisfies that $\rank N_1^f=m\leq 11$ and that 
$\nu^c_f(x)<2n-2m$ at any $x\in M^{2n}$. 
\vspace{1ex}

Since $\dim\Sal(\gamma(x))\leq 2m$, then
$\nu^c_f(x)<2n-\dim\Sal(\gamma(x))$ at any $x\in M^{2n}$.
Then by Theorem \ref{alglemma} the vector subbundle 
$\Omega\subset N_1^f$ has rank $2r\neq 0$ and 
$D_0^{d_0}=\mathcal{N}_c(\a_P)$ satisfies
\be\label{Gamma}
d_0\geq 2n-2m+4r.
\ee

We have that  $L_0^{2\ell_0}\neq 0$. If otherwise and 
since $D_0$ is $J$-invariant, then $D_0\subset\mathcal{N}_c(\a(x))$. 
Hence \eqref{Gamma} yields $\nu^c_f(x)>2n-2m$, which is a 
contradiction. It follows from \eqref{Gamma} that
\be\label{Gamma2}
d_0\geq 2n-2m+4\ell_0.
\ee

\begin{proposition}\po\label{step2} 
We have that $d_1\geq d_0-2(\ell_0-\ell_1)$ and if 
$\ell_1<\ell_0$ then $\ell_0-\ell_1\geq 2$.
\end{proposition}

\proof We assume that $\ell_0>\ell_1$ since equality yields
$D_1=D_0$, and thus the result holds trivially. From the 
definitions of $L_0$ and that of 
$R_1$ by the orthogonal decomposition $L_0=L_1\oplus R_1$, 
it follows that the bilinear form 
$\theta=\a_{R_1}|_{D_0\times TM}\colon D_0\times TM\to R_1$
satisfies $R_1=\Sal(\theta)$. From \eqref{condd1} we obtain that
\be\label{thetajinv}
\J_0\theta(S,X)=\theta(S,JX)=\theta(JS,X)
\ee
for any $S\in\Gamma(D_0)$ and $X\in\mathfrak{X}(M)$. 
From \eqref{thetajinv} it follows that
\be\label{D1}
D_1=\{S\in D_0\colon\theta(S,X)=0
\;\mbox{for all}\;X\in\mathfrak{X}(M)\}.
\ee

For $Z\in RE(\theta)$ we have by \eqref{thetajinv} that
$V_Z=\theta(Z,TM)$ is $\J_0$-invariant. Therefore
$\dim V_Z=s>0$ is even. We argue that $V_Z\neq R_1$. Since
$\mathcal{K}(X)$ is skew-symmetric then
$\mathcal{K}(X)R_1\subset R_1$ for any $X\in\mathfrak{X}(M)$. 
If $V_Z=R_1$, we have from part $(iv)$ of Lemma \ref{charK}  
that $\mathcal{K}(X)R_1\subset L_1$. 
But then $r_1=0$ where $\rank R_1=2r_1$, and this is a 
contradiction. It follows that $2r_1>s>0$ and, in particular, 
we already proved that $r_1=\ell_0-\ell_1\geq 2$.
 
We show next that it can only happen that 
$$
s=2\;\;\mbox{for}\;\;2\leq r_1\leq 3\;\;\mbox{and}\;\;
s=2,4\;\;\mbox{for}\;\;4\leq r_1\leq 5.
$$
If otherwise, we may have $s=4$ if $r_1=3$, $s=6$ if $r_1=4$ 
or $s=6,8$ if $r_1=5$. Let $k_0\geq 2$ be the minimal number 
of vectors $Z_1,\ldots,Z_{k_0}\in RE(\theta)$ such that
$R_1=\sum_{j=1}^{k_0}V_{Z_j}$. The dimension of $V_{Z_i}\cap V_{Z_j}$ 
is even.  Thus $k_0=2$ for $r_1=3,4$ and $k_0=2,3$ for $r_1=5$.  
It is easy to see that $U=\cap_{j=1}^{k_0}V_{Z_j}\neq 0$. 
On the other hand, part $(iv)$ of Lemma \ref{charK} gives that 
$\mathcal{K}(X)U\subset\cap_{j=1}^{k_0}V_{Z_j}^\perp\cap R_1=0$
for any $X\in\mathfrak{X}(M)$, and this is a contradiction.

Suppose that $Y\in\mathfrak{X}(M)$ satisfies $R_1=\theta(D_0,Y)$. 
Part $(iii)$ of Lemma \ref{charK} yields
\be\label{seqeq}
\mathcal{K}(X)\theta(S,Y)=\mathcal{K}(X)\a(S,Y)
=\mathcal{K}(Y)\a(S,X)=\mathcal{K}(Y)\theta(S,X)
\ee
for any $X,Y\in\mathfrak{X}(M)$ and $S\in\Gamma(D_0)$.
Then \eqref{ldspanR} and \eqref{seqeq} give
\be\label{supfirst}
\mathcal{K}(Y)R_1=R_1.
\ee
Now assume that there is $Y\in\mathfrak{X}(M)$ such that 
\eqref{supfirst} holds. From \eqref{seqeq} if $\theta(S,Y)=0$ 
then $\theta(S,X)=0$ for any $X\in\mathfrak{X}(M)$. 
It follows from \eqref{D1} that
$D_1=\{S\in D_0\colon \theta(S,Y)=0\}$,
and thus $d_0\leq d_1+2r_1=d_1+2(\ell_0-\ell_1)$ as  
required.
\vspace{1ex}

\noindent Case $s=2$. If $Z\in RE(\theta)$ let 
$\{X_i,JX_i\}_{1\leq i\leq n}$ be a basis  such that
$$
V_Z=\spa\{\theta(Z,X_1),\theta(Z,JX_1)\}
$$
and $\theta(Z,X_j)=\theta(Z,JX_j)=0$ if $j\geq 2$.
Given any $X=\sum_{i=1}^n a_iX_i+b_iJX_i$ and $S\in D_0$,
then by \eqref{firstst} and \eqref{thetajinv}
there is $T\in D_0$ such that
$$
\theta(S,X)=\sum_{i=1}^n a_i\theta(S,X_i)
+\sum_{i=1}^n b_i\theta(S,JX_i)
=\theta(T,X_1).
$$
Since $R_1=\Sal(\theta)$ then $R_1=\theta(D_0,X_1)$, and the 
result follows as argued above.
\vspace{1ex}

\noindent Case $s=4$. If $Z\in RE(\theta)$
there is a basis $\{X_i,JX_i\}_{1\leq i\leq n}$ such that
$$
V_Z=\spa\{\theta(Z,X_1),\theta(Z,JX_1),
\theta(Z,X_2),\theta(Z,JX_2)\}
$$
and $\theta(Z,X_j)=\theta(Z,JX_j)=0$ if $j\geq 3$.
Given any $X=\sum_{i=1}^n a_iX_i+b_iJX_i$ and $S\in D_0$
we have from \eqref{firstst} and \eqref{thetajinv} that
there are vectors $T_1,T_2\in D_0$ such that
$$
\theta(S,X)=\sum_{i=1}^n a_i\theta(S,X_i)
+\sum_{i=1}^n b_i\theta(S,JX_i)
=\theta(T_1,X_1)+\theta(T_2,X_2).
$$
Since $R_1=\Sal(\theta)$ then
\be\label{imatheta}
R_1=\theta(D_0,X_1)+\theta(D_0,X_2). 
\ee
As seen above, the only cases that have to be considered
is when $\theta(D_0,X_j)\neq R_1$ for $j=1,2$ and 
$\mathcal{K}(Y)R_1\neq R_1$ for any 
$Y\in\mathfrak{X}(M)$. From \eqref{thetajinv}
we have that $\theta(D_0,X_j)$ es $\J_0$-invariant.
Hence $N_j=\{S\in D_0\colon \theta(S,X_j)=0\}$ satisfies
\be\label{estimnj}
\dim N_j\geq d_0-2r_1+2,\;j=1,2.
\ee

By part $(iii)$ of Lemma \ref{charK} we have 
$\mathcal{K}(X)\theta(Z,X_j)=\mathcal{K}(X_j)\theta(Z,X)$, 
$j=1,2$, for any $X\in\mathfrak{X}(M)$ and $Z\in D_0$.
Then \eqref{ldspanR} and \eqref{imatheta} yield
\be\label{step2r1}
R_1=\Ima\mathcal{K}(X_1)+\Ima\mathcal{K}(X_2).
\ee
Let $\mathcal{K}\colon\mathfrak{X}(M)\times\Gamma(R_1)\to\Gamma(R_1)$
be the bilinear form defined by $\mathcal{K}(X,\eta)=\mathcal{K}(X)\eta$. 
We may assume that $X_1,X_2\in RE(\mathcal{K})$, and hence
$\dim\Ima\mathcal{K}(X_1)=\dim\Ima\mathcal{K}(X_2)$.

We claim that $\dim\Ima\mathcal{K}(X_j)=2r_1-2$ for $j=1,2$.
We only argue for the most difficult case $r_1=5$. Then $\ell_0=5$
and $\ell_1=0$. If the claim does not hold, it follows  
from \eqref{step2r1} that
$\dim\Ima\mathcal{K}(X_j)=6$, $j=1,2$. Hence
$\dim(\Ima\mathcal{K}(X_1)\cap\Ima\mathcal{K}(X_2))=2$.
Since $\mathcal{K}(X)$ is skew-symmetric then
$(\Ima\mathcal{K}(X))^\perp=\ker\mathcal{K}(X)$, and hence
\be\label{step2inter}
\ker\mathcal{K}(X_1)\cap\ker\mathcal{K}(X_2)=0.
\ee
If $aX_1+bX_2\in RE(\mathcal{K})$ with $a\neq 0\neq b$, then
$V=\ker\mathcal{K}(aX_1+bX_2)$ satisfies $\dim V=4$.
We have that $\mathcal{K}(X_1)V=\mathcal{K}(X_2)V$. Hence 
$\mathcal{K}(X_1)V
\subset\Ima\mathcal{K}(X_1)\cap\Ima\mathcal{K}(X_2)$.
Then $\ker\mathcal{K}(X_1)\cap V\neq 0$
in contradiction with \eqref{step2inter}.
\vspace{1ex}

We have
$$
\mathcal{K}(X_1)\theta(N_1,X_2)=\mathcal{K}(X_1)\a(N_1,X_2)
=\mathcal{K}(X_2)\a(N_1,X_1)=\mathcal{K}(X_2)\theta(N_1,X_1)=0.
$$
Since $\dim\Ima\mathcal{K}(X_j)=8$, $j=1,2$, then
$\dim\Ima\theta(N_1,X_2)\leq\dim\ker\mathcal{K}(X_1)=2$.
We obtain from \eqref{estimnj} that $N=N_1\cap N_2$ satisfies 
$\dim N\geq d_0-10$, and from part $(iii)$ of Lemma \ref{charK} 
that
$$
\mathcal{K}(X_1)\theta(S,X)=\mathcal{K}(X_2)\theta(S,X)=0
$$
for any $S\in N$ and $X\in\mathfrak{X}(M)$. Then \eqref{step2inter}
gives $\theta(S,X)=0$ and \eqref{D1} that $N\subset D_1$.\qed

\begin{proposition}\po\label{step3}
We have that $\ell\geq 1$ and $d\geq d_0-4(\ell_0-\ell)$.
\end{proposition}

\proof We may assume that $\ell_0>\ell$. If $L_0=L\oplus R$ is
an orthogonal decomposition, we obtain from the definition of $L_0$ 
that $\bar\theta=\a_R|_{TM\times D_0}\colon TM\times D_0\to R$
satisfies $R=\Sal(\bar\theta)$.
Set $\upsilon=\dim\bar\theta_X(D_0)$ for $X\in RE(\bar\theta)$.  
From \eqref{condd00} we obtain
$$
\J_0\a_R(X,S)=\a_R(X,JS)
$$
for any $X\in\mathfrak{X}(M)$ and $S\in\Gamma(D_0)$. Thus 
$\bar\theta_X(D_0)$ is $\J_0$-invariant and $\upsilon$ is even.

Let the vector subspace $V^k\subset TM$ be of minimum dimension    
so that $R=\Sal(\bar\theta|_{V^k\times D_0})$.
If $X_1,\ldots, X_k\in RE(\bar\theta)$ span $V^k$ then 
$R=\sum_{i=1}^k \bar\theta_{X_i}(D_0)$. 
Since $\bar\theta_{X_i}(D_0)\cap\bar\theta_{X_j}(D_0)$ is 
$\J_0$-invariant, thus of even dimension, then  
$2(\ell_0-\ell)\geq\upsilon+2(k-1)$. Hence
\be\label{rhoineq}
2\leq\upsilon\leq 2(\ell_0-\ell-k+1).
\ee 

From \eqref{firstst} we obtain that
$\Sal(\bar\theta|_{TM\times\ker\bar\theta_{X_j}})
\subset\bar\theta_{X_j}(D_0)$. By the definition of $k$, we have 
$$
\bar\theta_{X_i}(\ker\bar\theta_{X_j})
\subsetneq\bar\theta_{X_j}(D_0)\;\;\mbox{if}\;\;i\neq j.
$$
Since $\dim\ker\bar\theta_{X_j}=d_0-\upsilon$, it follows that 
$$
\dim(\ker\bar\theta_{X_i}\cap\ker\bar\theta_{X_j})
\geq d_0-2\upsilon+2\;\;\mbox{if}\;\;i\neq j.
$$
If $k\geq 3$, we have 
$$
\bar\theta_{X_i}(\ker\bar\theta_{X_j}\cap\ker\bar\theta_{X_r})
\subsetneq\bar\theta_{X_j}(D_0)
\;\;\mbox{if}\;\;i\neq j\neq r\neq i.
$$
Hence
$$
\dim(\ker\bar\theta_{X_i}\cap\ker\bar\theta_{X_j}\cap\ker\bar\theta_{X_r})
\geq \dim(\ker\bar\theta_{X_i}\cap\ker\bar\theta_{X_j})-\upsilon+2
\geq d_0-3(\upsilon-2)-2.
$$
Reiterating the argument gives that
$G=\mathcal{N}(\bar\theta|_{V^k\times D_0})
=\cap_{i=1}^k\ker\bar\theta_{X_i}$ satisfies 
\be\label{dimg0}
\dim G\geq d_0-k(\upsilon-2)-2.
\ee

We claim that 
\be\label{defdalt}
D=G\cap D_1.
\ee 
The definition of $D$ yields
$0=\a_{L^\perp}(X,S)=\a_R(X,S)+\a_{L_0^\perp}(X,S)$
for any $X\in V^k$ and $S\in\Gamma(D)$. Thus $S\in G$, 
which gives one inclusion.  

Taking the  $L_0^\perp$-component of 
$(\nabla_Z^\perp\a)(X,Y)=(\nabla_Y^\perp\a)(X,Z)$ 
and using the Lemma \ref{liebracketprop} we obtain
$$
(\nabla^\perp_Z\a(X,Y))_{L_0^\perp}
=(\nabla^\perp_Y\a(X,Z))_{L_0^\perp}
$$ 
for any $Z\in\Gamma(D_1)$, $X\in\mathfrak{X}(M)$  
and $Y\in\Gamma(D_0)$. The definition of $L$ gives 
$(\nabla_Z^\perp\zeta)_{L_0^\perp}=0$ for any 
$Z\in\Gamma(D_0)$ and $\zeta\in\Gamma(L)$. Since 
$\a(X,Y)\in\Gamma(L_0)$ and 
$\a(X,Z)\in\Gamma(L_1)\subset\Gamma(L_0)$, then
$$
(\nabla^\perp_Z\bar\theta(X,Y))_{L_0^\perp}
=(\nabla^\perp_Z\a_R(X,Y))_{L_0^\perp}
=(\nabla^\perp_Y\a_R(X,Z))_{L_0^\perp}
=(\nabla^\perp_Y\bar\theta(X,Z))_{L_0^\perp}=0
$$
for any $Z\in\Gamma(G\cap D_1)$, $X\in\Gamma(V^k)$ and 
$Y\in\Gamma(D_0)$. Since $R=\Sal(\bar\theta|_{V^k\times D_0})$ then
\be\label{parallelstep3}
(\nabla_Z^\perp\xi)_{L_0^\perp}=0
\ee 
for any $Z\in\Gamma(G\cap D_1)$ and $\xi\in\Gamma(R)$.
Then the  above Codazzi equation  and \eqref{parallelstep3} yield
$$
(\nabla_Y^\perp\a_R(X,Z))_{L_0^\perp}
=(\nabla_Z^\perp\a_R(X,Y))_{L_0^\perp}=0
$$
for any $X\in\mathfrak{X}(M)$, $Y\in \Gamma(D_0)$ and 
$Z\in\Gamma(G\cap D_1)$. Hence 
$$
\nabla_Y^\perp\a_R(X,Z)=
\nabla_Y^\perp\a_{L^\perp\cap L_1}(X,Z)
\in\Gamma(L_0)
$$
for any $X\in\mathfrak{X}(M)$, $Y\in \Gamma(D_0)$ and 
$Z\in\Gamma(G\cap D_1)$.  
Since $G$ is $\J_0$-invariant, then also $G\cap D_1$ is 
$\J_1$-invariant and thus
$$
\J_1\a_{L^\perp\cap L_1}(X,Z)
=\a_{L^\perp\cap L_1}(X,JZ)\in\Gamma(L_0)
$$
for any $X\in\mathfrak{X}(M)$ and $Z\in\Gamma(G\cap D_1)$. 
Hence also
$\nabla_Y^\perp\J_1\a_{L^\perp\cap L_1}(X,Z)\in\Gamma(L_0)$.
Then \eqref{defl2} gives that $\a_{L^\perp\cap L_1}(X,Z)\in\Gamma(L)$,
that is, that $\a_{L^\perp\cap L_1}(X,Z)=0$ 
for any $X\in\mathfrak{X}(M)$ and $Z\in\Gamma(G\cap D_1)$. Since
$$
\a_{L^\perp}(X,Z)=\a_{L^\perp\cap L_1}(X,Z)+\a_{L_1^\perp}(X,Z)=0
$$
for any $X\in\mathfrak{X}(M)$ and $Z\in\Gamma(G\cap D_1)$, then 
$Z\in\Gamma(D)$ thus proving the claim.  
\vspace{1ex}

Since $G\subset D_0$ we have from \eqref{defdalt} that
$d_0\geq\dim G+d_1-d$. Then \eqref{dimg0} gives
\be\label{dimg1}
d\geq d_1-k(\upsilon-2)-2.
\ee 

We have that $\hat{\theta}=\a_{L^\perp\cap L_1}|_{V^k\times D_1}
\colon V^k\times D_1\to L^\perp\cap L_1$ satisfies
\be\label{hattheta}
D=\mathcal{N}(\hat\theta).
\ee
In fact, that $D\subset\mathcal{N}(\hat\theta)$ follows 
from the definition of $D$. Take $S\in\mathcal{N}(\hat\theta)\subset D_1$. 
By \eqref{defdalt} it suffices to show that $S\in G$, which follows from
$$
\a_R(X,S)=\a_{L^\perp\cap L_1}(X,S)
+\a_{L_0\cap L_1^\perp}(X,S)=\hat\theta(X,S)=0
$$
for any $X\in V^k$.

Denote $\upsilon_1=\dim\hat\theta_Y(D_1)\leq 2(\ell_1-\ell)$
when $Y\in RE(\hat\theta)$. If $Y_1,\ldots,Y_k\in RE(\hat\theta)$ 
span $V^k$ then \eqref{firstst} gives 
$\hat\theta_{Y_i}(\ker\hat\theta_{Y_j})\subset \hat\theta_{Y_j}(D_1)$
if $i\neq j$. Thus
$$
d_1-\upsilon_1=\dim\ker\hat\theta_{Y_j}\leq
\dim(\ker\hat\theta_{Y_i}\cap\ker\hat\theta_{Y_j})
+\upsilon_1\;\;\mbox{if}\;\;i\neq j.
$$
If $k\geq 3$, we have
$$
\hat\theta_{Y_i}(\ker\hat\theta_{Y_j}\cap\ker\hat\theta_{Y_r})
\subset\hat\theta_{Y_j}(D_1)
\;\;\mbox{if}\;\;i\neq j\neq r\neq i.
$$
Thus
$$
d_1-2\upsilon_1\leq
\dim(\ker\hat\theta_{Y_i}\cap\ker\hat\theta_{Y_j})-\upsilon_1 
\leq \dim(\ker\hat\theta_{Y_i}\cap\ker\hat\theta_{Y_j}\cap\ker\hat\theta_{Y_r}).
$$
Reiterating the argument gives
$$
d_1-k\upsilon_1\leq\dim\cap_{i=1}^k\ker\hat\theta_{Y_i}.
$$
Then by \eqref{hattheta} we have
$$
d=\dim\cap_{i=1}^k\ker\hat\theta_{Y_i}\geq
d_1-k\upsilon_1\geq d_1-2k(\ell_1-\ell).
$$ 
This and \eqref{dimg1} yield
$$
d\geq d_1-\min\{k(\upsilon-2)+2,2k(\ell_1-\ell)\}.
$$
Now Proposition \ref{step2} gives that
$$
d\geq d_0-2(\ell_0-\ell_1)-\min\{k(\upsilon-2)+2,2k(\ell_1-\ell)\}.
$$
Therefore, to conclude the proof it suffices to show that
$$
\min\{k(\upsilon-2)+2,2k(\ell_1-\ell)\}\leq 2\ell_0+2\ell_1-4\ell.  
$$
If this does not hold and since $\upsilon$ is even, then
$$
k(\upsilon-2)\geq 2\ell_0+2\ell_1-4\ell
\;\;\mbox{and}\;\;
(2k-4)(\ell_1-\ell)\geq 2\ell_0-2\ell_1+2.
$$
Since $\ell_0>\ell$, then the first inequality 
gives $\upsilon\geq 4$ whereas the second that $k\geq 3$. 
On the other hand, since $2\ell_0-2\ell\leq 2\ell_0\leq 10$ 
from \eqref{rhoineq} it remains to analyze the cases in which
$(k,\upsilon)$ is either $(3,4)$, $(4,4)$ and $(3,6)$. For the 
last two cases we have from \eqref{rhoineq} that 
$2(\ell_0-\ell)=10$ and thus $\ell_0=5$ and $\ell=0$.  
For $(4,4)$ the first inequality gives a contradiction.  
For $(3,6)$ the first inequality gives $0\leq\ell_1\leq 1$, 
and then the second yields a contradiction. For the case $(3,4)$ 
we have from \eqref{rhoineq} that $8\leq 2(\ell_0-\ell)\leq 10$,  
and the first inequality gives a contradiction.

Finally, if $\ell=0$ it follows from \eqref{Gamma2} 
that $\nu^c_f=d\geq d_0-4\ell_0\geq 2n-2m$, and this is a 
contradiction.\qed

\begin{proposition}\po\label{destimate}
If $\ell_1=\ell_0\leq\ell + 2$ then $d\geq d_0-2(\ell_0-\ell)$.
\end{proposition}

\proof We assume that $\ell_0>\ell$ since the 
result holds trivially otherwise. We argue that
\be\label{shwoed}
L^\perp\cap L_0
=\spa\{(\nabla_S^\perp\eta)_{L_0},\J_0(\nabla_S^\perp\eta)_{L_0}
\colon S\in\Gamma(D_0)\;\mbox{and}\;\eta\in\Gamma(L_0^\perp)\}.
\ee
We have that
\be\label{abeq}
\<(\nabla_S^\perp\eta)_{L_0},\xi\>=-\<\eta,\nabla_S^\perp\xi\>
\;\;\mbox{and}\;\;
\<\J_0(\nabla_S^\perp\eta)_{L_0},\xi\>
=\<\eta,\nabla_S^\perp\J_0\xi\>
\ee
for any $S\in\Gamma(D_0)$, $\eta\in\Gamma(L_0^\perp)$ and 
$\xi\in\Gamma(L_0)$.  The left-hand-side of \eqref{abeq} vanishes 
if $\xi\in\Gamma(L)$. Hence
$(\nabla_S^\perp\eta)_{L_0},\J_0(\nabla_S^\perp\eta)_{L_0}
\in\Gamma(L^\perp\cap L_0)$, which gives one inclusion. 

Now assume that $\xi\in\Gamma(L^\perp\cap L_0)$ satisfies
$$
\<(\nabla_S^\perp\eta)_{L_0},\xi\>
=\<\J_0(\nabla_S^\perp\eta)_{L_0},\xi\>=0
$$
for any $S\in\Gamma(D_0)$ and $\eta\in\Gamma(L_0^\perp)$. 
Then \eqref{abeq} yields that 
$\nabla_S^\perp\xi,\nabla_S^\perp\J_0\xi\in\Gamma(L_0)$
for any $S\in\Gamma(D_0)$. Hence $\xi\in\Gamma(L)$ and 
thus $\xi=0$, proving the other inclusion.

The tensor 
$\psi\colon\Gamma(L_0^\perp)\times\Gamma(D_0)\to\Gamma(L^\perp\cap L_0)$
defined by $\psi(\eta,S)=\psi_\eta S=(\nabla_S^\perp\eta)_{L_0}$
satisfies $\Sal(\psi)+\J_0\Sal(\psi)=L^\perp\cap L_0$. 
By Lemma \ref{liebracketprop} the  distribution $D_0$ is integrable, 
and hence the Codazzi equation
$(\nabla_S A)(\eta,T)=(\nabla_T A)(\eta,S)$ gives that
$A_{\psi_\eta S}T=A_{\psi_\eta T}S$
for any $S,T\in\Gamma(D_0)$ and $\eta\in\Gamma(L_0^\perp)$.  
Then, we obtain using \eqref{condd00} that
$$
JA_{\psi_\eta S}T=A_{\J_0\psi_\eta T}S
=-A_{\psi_\eta T}JS=JA_{\psi_\eta T}S.
$$
Hence
\be\label{kern}
\ker\psi_\eta+J\ker\psi_\eta
\subset\ker A_{\Ima\psi\eta}\cap\ker A_{\J_0\Ima\psi\eta}.
\ee

Assume that there is $\eta\in\Gamma(L_0^\perp)$ such that
$\Ima\psi_\eta+\J_0\Ima\psi_\eta=L^\perp\cap L_0$.
By \eqref{shwoed} this is necessarily the case if 
$\ell_0-\ell=1$. From \eqref{kern} we have 
$\ker\psi_\eta\subset D$, and it follows that
$d\geq\dim\ker\psi_\eta\geq d_0-2(\ell_0-\ell)$.

The remaining case to be considered is when $\ell_0-\ell=2$ 
and there are $\eta_j\in\Gamma(L_0^\perp)$, $j=1,2$, with 
$\dim\psi_{\eta_j}(D_0)\leq 2$ such that
$$
\psi_{\eta_1}(D_0)+\J_0\psi_{\eta_1}(D_0)+
\psi_{\eta_2}(D_0)+\J_0\psi_{\eta_2}(D_0)=L^\perp\cap L_0.
$$
But then 
$d\geq\dim(\ker\psi_{\eta_1}\cap\ker\psi_{\eta_2})\geq d_0-4$.
\qed

\section{The proofs}

The proof of both theorems presented in the Introduction 
ensues by amalgamating the three established theorems in 
this section.
It is noteworthy that the assertions of the latter set 
of theorems exhibit a slightly higher degree of generality, 
for instance, the Kaehler extensions are explicitly constructed 
through the application of Theorem \ref{develop}.
\vspace{2ex}

We recall that all results 
hold on connected components of the open dense subset of the 
submanifold where the various vector subspaces that has been 
under consideration maintain constant dimension, thus forming 
vector bundles.
\vspace{2ex}

Let $\phi\colon\mathfrak{X}(M)\times\Gamma(f_*TM\oplus L_1)
\to\Gamma(L_0^\perp\oplus L_0^\perp)$ be the bilinear 
tensor given~by
$$
\phi(Y,\zeta)=((\tilde\nabla_Y\zeta)_{L_0^\perp},
(\tilde\nabla_Y\hat\J_1\zeta)_{L_0^\perp}),
$$ 
where $\hat{\J}_1=\hat{\J}_0|_{f_*TM\oplus L_1}$ and 
$\hat\J_0\in\Gamma(\text{Aut}(f_*TM\oplus L_0))$ is defined by 
$$
\hat{\J}_0(f_*X+\eta)=f_*JX+\J_0\eta.
$$
Thus
$$
\phi(Y,f_*X+\xi)=(\a_{L_0^\perp}(Y,X)+(\nap_Y\xi)_{L_0^\perp},
\a_{L_0^\perp}(Y,JX)+(\nap_Y\hat\J_1\xi)_{L_0^\perp}).
$$
Moreover, since $\mathcal{N}(\phi)$ is $\hat\J_1$-invariant 
then $\nu(\phi)$ is even.  

\begin{proposition}\po\label{tensorphi}
The bilinear form $\phi(x)$  is flat at any 
$x\in M^{2n}$ with respect to the inner product 
induced on 
$L_0^\perp\oplus L_0^\perp\subset N_fM\oplus N_fM$. 
\end{proposition}

\proof It follows from \eqref{condd00} and \eqref{parallell0} that
\be\label{hatj0}
\hat\J_0((\tilde\nabla_X\zeta)_{f_*TM\oplus L_0})
=(\tilde\nabla_X\hat\J_0\zeta)_{f_*TM\oplus L_0}
\ee
for any $X\in\mathfrak{X}(M)$ and $\zeta\in\Gamma(f_*TM\oplus L_1)$.
Denoting by $\tilde{R}$ the curvature tensor of the ambient space,
we have from $\<\tilde{R}(X,Y)\zeta,\varsigma\>=0$ that
$$
\<\tilde\nabla_Y\zeta,\tilde\nabla_X\varsigma\>
-\<\tilde\nabla_X\zeta,\tilde\nabla_Y\varsigma\>
=X\<\tilde\nabla_Y\zeta,\varsigma\>-Y\<\tilde\nabla_X\zeta,\varsigma\>
-\<\tilde\nabla_{[X,Y]}\zeta,\varsigma\>,
$$
whereas from $\<\tilde{R}(X,Y)\hat\J_0\zeta,\hat\J_0\varsigma\>=0$ that
\begin{align*}
\<\tilde\nabla_Y\hat\J_0\zeta,\tilde\nabla_X\hat\J_0\varsigma\>
-\<\tilde\nabla_X\hat\J_0\zeta,\tilde\nabla_Y\hat\J_0\varsigma\>&=
X\<\tilde\nabla_Y\hat\J_0\zeta,\hat\J_0\varsigma\>
-Y\<\tilde\nabla_X\hat\J_0\zeta,\hat\J_0\varsigma\>\\
&\quad-\<\tilde\nabla_{[X,Y]}\hat\J_0\zeta,\hat\J_0\varsigma\>
\end{align*}
for any $X,Y\in\mathfrak{X}(M)$ and 
$\zeta,\varsigma\in\Gamma(f_*TM\oplus L_1)$.
Since by \eqref{hatj0} the right-hand-sides of the two preceding
equations coincide, then
$$
\<\tilde\nabla_Y\zeta,\tilde\nabla_X\varsigma\>
-\<\tilde\nabla_X\zeta,\tilde\nabla_Y\varsigma\>
=\<\tilde\nabla_Y\hat\J_0\zeta,\tilde\nabla_X\hat\J_0\varsigma\>
-\<\tilde\nabla_X\hat\J_0\zeta,\tilde\nabla_Y\hat\J_0\varsigma\>
$$
for any $X,Y\in\mathfrak{X}(M)$ and
$\zeta,\varsigma\in\Gamma(f_*TM\oplus L_1)$. 
By \eqref{hatj0} we also have that
$$
\<(\tilde\nabla_Y\zeta)_{f_*TM\oplus L_0},
(\tilde\nabla_X\varsigma)_{f_*TM\oplus L_0}\>
=\<(\tilde\nabla_Y\hat\J_0\zeta)_{f_*TM\oplus L_0},
(\tilde\nabla_X\hat\J_0\varsigma)_{f_*TM\oplus L_0}\>.
$$
It follows that
$$
\lp\phi(Y,\zeta),\phi(X,\varsigma)\rp
-\lp\phi(X,\zeta),\phi(Y,\varsigma)\rp=0
$$
for any  $X,Y\in\mathfrak{X}(M)$ and 
$\zeta,\varsigma\in\Gamma(f_*TM\oplus L_1)$.\qed 

\begin{proposition}\po\label{phiestimates}
Let $f\colon M^{2n}\to\R^{2n+p}$, $p\leq n-1$ and $p\leq 11$, 
be a real Kaehler submanifold satisfying $\nu^c_f(x)<2n-2p$ 
at any $x\in M^{2n}$. If $p\leq 2\ell_0+5$ then at any 
$x\in M^{2n}$ at least one of the following inequalities 
holds:
\be\label{phiest1}
\nu(\phi)\geq 2n-2p+4\ell_0+2\ell_1
\ee
and
\be\label{phiest2}
d_0\geq 2n-2p+4\ell_0+4.
\ee
\end{proposition}

\proof Since $\phi_X(\delta)=(\xi,\bar\xi)$ if and only if
$\phi_X(\hat\J_1\delta)=(\bar\xi,-\xi)$, then the vector bundles 
$\Ima\phi_X$ and $\Sal(\phi)$ are either zero or of even rank.  
Set $\U(X)=\Ima\phi_X\cap(\Ima\phi_X)^\perp$ where 
$X\in\mathfrak{X}(M)$ and $\U=\Sal(\phi)\cap\Sal(\phi)^\perp
\subset L_0^\perp\oplus L_0^\perp$. We have that $\U(X)$ and $\U$  
are zero or of even rank. In fact, this it is easily seen since 
\be\label{argue3}
(\xi,\bar\xi)\in\U(X)\;\;\mbox{if and only if}\;\;
(\bar\xi,-\xi)\in\U(X).
\ee

We have that $L_0^{2\ell_0}\subset\Omega^{2r}$ where 
$\Omega^{2r}=\pi_1(\U)$. 
If $\ell_0<r$ 
then \eqref{Gamma} gives \eqref{phiest2}. 
Therefore, in the sequel we assume that $L_0=\Omega$.

We show that \eqref{phiest2} also holds if $\U\neq 0$.  
In fact, let us consider the decomposition
$L_0^\perp\oplus L_0^\perp=\U\oplus\hat\U\oplus\mathcal{V}_0$,
where $\Sal(\phi)\subset\U\oplus\mathcal{V}_0$ and the vector 
subspace $\mathcal{V}_0=(\U\oplus\hat\U)^\perp$ is nondegenerate.  
Notice that we have
$\Sal(\gamma_{L_0^\perp})=\Sal(\gamma_{\Omega^\perp})\subset\Sal(\phi)$.
Since Theorem \ref{alglemma} gives that the vector subspace 
$\Sal(\gamma_{L_0^\perp})$ is nondegenerate, then the map 
$\pi_{\mathcal{V}_0}|_{\Sal(\gamma_{L_0^\perp})}$ is injective, where 
$\pi_{\mathcal{V}_0}\colon L_0^\perp\oplus L_0^\perp\to\mathcal{V}_0$ 
denotes taking the $\mathcal{V}_0$-component. Being $\dim\U=\dim\hat\U$ 
even then 
$\dim\Sal(\gamma_{L_0^\perp})\leq\dim\mathcal{V}_0\leq 2p-4\ell_0-4$, 
and we obtain \eqref{phiest2} from \eqref{condkerd0} and 
Theorem~\ref{alglemma}. Therefore, in the sequel we also 
assume that $\U=0$.

If $\phi=0$ which, in particular, happens when $L_0=N_fM$, then 
\eqref{phiest1} holds trivially. Moreover, since $\nu^c_f<2n-2p$ 
we have from \eqref{Gamma} that $L_0\neq 0$. Therefore, in 
the sequel we also assume that $\phi\neq 0$ and $0<2\ell_0<p$.

By Proposition \ref{facts},  the subset
$RE^\#(\phi)\subset RE(\phi)$ of $T_xM$ defined by
$$
RE^\#(\phi)=\{X\in RE(\phi)\colon \dim\U(X)=\tau\},
$$
where $\tau=\min\{\dim\U(X)\colon X\in RE(\phi)\}$,
is open and dense.
Fix $X\in RE^\#(\phi)$ and denote $\kappa=\dim\Ima\phi_X$.
Then either $\tau=0$ or $\tau\leq \kappa$ are both even.
By Sublemma $2.3$ in \cite{CD} or Corollary $4.3$ 
in \cite{DT} there is the decomposition
\be\label{cod6decomp}
L_0^\perp\oplus L_0^\perp=
\U^\tau(X)\oplus\hat\U^\tau(X)\oplus
\mathcal{V}^{p-2\ell_0-\tau,p-2\ell_0-\tau}
\ee
such that $\Ima\phi_X\subset\U(X)\oplus
\mathcal{V}$, $\mathcal{V}
=(\U(X)\oplus\hat\U(X))^\perp$ 
is a nondegenerate vector subspace and   
$0\leq\tau\leq \kappa\leq 2p-4\ell_0-\tau$.
In fact, we cannot have that $\tau=\kappa$. 
If otherwise, the vector subspace $\Ima\phi_X$ 
is isotropic and thus
$$
\lp\phi(X,\delta),\phi(X,\xi)\rp=0
$$
for any $X\in RE^\#(\phi)$ and 
$\delta,\xi\in\Gamma(f_*TM\oplus L_1)$. Then the same 
holds for any $X\in\mathfrak{X}(M)$, and being $\phi$ 
flat by Proposition \ref{tensorphi}, we have that
$$
2\lp\phi(Y,\delta),\phi(Z,\xi)\rp=
\lp\phi(Y+Z,\delta),\phi(Y+Z,\xi)\rp=0
$$
for any $Y,Z\in\mathfrak{X}(M)$ and
$\delta,\xi\in\Gamma(f_*TM\oplus L_1)$. Since 
$\phi\neq 0$ and $\U=0$, this is a contradiction. 

If $\tau=0$ then \eqref{secondst} gives that
$\mathcal{N}(\phi)=\ker\phi_X$, and hence
$$
\nu(\phi)=\dim\ker\phi_X\geq 2n-2p+4\ell_0+2\ell_1
$$
which is \eqref{phiest1}. Therefore, since by assumption
$p\leq 2\ell_0+5$, it remains to consider the case 
\be\label{phikerest2}
0<\tau<\kappa\leq 2p-4\ell_0-\tau\leq 10-\tau.
\ee

Let $\hat\phi$ denote the $\hat\U(X)$-component of $\phi$.
Since $\U=0$ thus $\Sal(\hat\phi)=\hat\U(X)$. 
Moreover, given $Y\in RE(\hat\phi)$ we have that
$\kappa_0=\dim\Ima\hat\phi_Y$ is even. In fact, by 
\eqref{argue3} there 
is a basis $\{(\xi_j,\bar{\xi}_j),(\bar{\xi}_j,-\xi_j),
1\leq j\leq m\}$  of $\U(X)$ where $\tau=2m$. 
From \eqref{cod6decomp} we have
$$
\begin{cases}
\lp\hat\phi(Y,\delta),(\xi_s,\bar\xi_s)\rp
=-\lp\hat\phi(Y,\hat\J_1\delta),(\bar\xi_s,-\xi_s)\rp
\vspace{1ex}\\
\lp\hat\phi(Y,\delta),(\bar\xi_s,-\xi_s)\rp
=\lp\hat\phi(Y,\hat\J_1\delta),(\xi_s,\bar\xi_s)\rp
\end{cases}
$$
for any $\delta\in\Gamma(f_*TM\oplus L_1)$ and $1\leq s\leq m$. 
Thus $T\in\text{Aut}(\Ima\hat\phi_Z)$ that is given by
$T\hat\phi_Z(\delta)=\hat\phi_Z(\hat\J_1\delta)$ is well 
defined and satisfies $T^2=-I$. Hence $\kappa_0$ is even.

Assume that $\kappa_0=\tau$, which is the case when 
$\tau=2$. Thus $\Ima\hat\phi_Y=\hat\U(X)$ for any 
$Y\in RE(\hat\phi)$. Then,
given $Y\in RE(\hat\phi)$ we have by \eqref{secondst}
and \eqref{cod6decomp} that
$$
\phi(Z,\eta)=0\;\;\mbox{if and only if}\;\;
\lp\phi(Z,\eta),\hat\phi(Y,\delta)\rp=0
$$
for any $Z\in\mathfrak{X}(M)$, $\eta\in\Gamma(\ker\phi_X)$ and
$\delta\in\Gamma(f_*TM\oplus L_1)$. 
Set $N=\ker\phi_Y|_{\ker\phi_X}$ where  
$\phi_Y|_{\ker\phi_X}\colon\ker\phi_X\to\U(X)$
is well defined by \eqref{secondst}. Then \eqref{secondst},
\eqref{cod6decomp} and flatness of $\phi$ yield
$$
\lp\phi(Z,\zeta),\hat\phi(Y,\delta)\rp
=\lp\phi(Z,\zeta),\phi(Y,\delta)\rp
=\lp\phi(Z,\delta),\phi(Y,\zeta)\rp=0
$$
for any $Z\in\mathfrak{X}(M)$, $\zeta\in\Gamma(N)$ and
$\delta\in\Gamma(f_*TM\oplus L_1)$.
Thus $N=\mathcal{N}(\phi)$. It follows using 
\eqref{phikerest2} that
$$
\nu(\phi)=\dim N
\geq\dim\ker\phi_X-\tau\geq 2n+2\ell_1-\kappa
-\tau\geq 2n-2p+4\ell_0+2\ell_1,
$$
which is \eqref{phiest1}.

By \eqref{phikerest2} it remains to argue for the case 
$\kappa_0=2$, $\tau=4$ and $\kappa=6$. Since the vector subspace 
$\Ima\hat\phi_{Y_1}\cap\Ima\hat\phi_{Y_2}$ is $T$-invariant,  
there are $Y_1,Y_2\in RE^\#(\phi)\cap RE(\hat\phi)$ such that
$$
\hat\U(X)=\Ima\hat\phi_{Y_1}\oplus\Ima\hat\phi_{Y_2}.
$$
Given $Z\in\mathfrak{X}(M)$ and $\eta\in\Gamma(\ker\phi_X)$,
we obtain from \eqref{secondst} that
$$
\phi(Z,\eta)=0
\;\;\mbox{if and only if}\;\;
\lp\phi(Z,\eta),\hat\phi(Y_i,\delta)\rp=0,\;\;i=1,2,
$$
for any $\delta\in\Gamma(f_*TM\oplus L_1)$. We have that
$\phi_{Y_i}|_{\ker\phi_X}\colon \ker\phi_X\to\U(X)$, $i=1,2$,
satisfies $\dim\phi_{Y_i}(\ker\phi_X)\leq 2$. If otherwise, 
say $\phi_{Y_1}(\ker\phi_X)=\U(X)$.  From
$$
\lp\hat\phi_{Y_1}\delta,\xi\rp=
\lp\phi_{Y_1}\delta,\xi\rp=0
$$
for $\delta\in\Gamma(f_*TM\oplus L_1)$ and $\xi\in\U(X)\cap\U(Y_1)$, 
we obtain $\dim\U(X)\cap\U(Y_1)\leq 2$. Consider the decomposition
$$
\U(Y_1)=\U(X)\cap\U(Y_1)\oplus\mathcal{V}_1.
$$
Fix $\zeta\in\mathcal{V}_1$. We have that $\lp\phi_{Y_1}\eta,\zeta\rp=0$
for any $\eta\in\Gamma(\ker\phi_X)$ because $\zeta\in\U(Y_1)$.
Since $\phi_{Y_1}(\ker\phi_X)=\U(X)$ 
then $\zeta=\zeta_1+\zeta_2$ from \eqref{cod6decomp}, where
$\zeta_1\in\U(X)$ and $\zeta_2\in\mathcal{V}^{1,1}$ with $\zeta_2\neq 0$.  
Thus the map $\pi_0\colon\mathcal{V}_1\to\mathcal{V}^{1,1}$ 
given by $\pi_0(\zeta)=\zeta_2$
is injective.  Moreover, we have that 
$\pi_0(\mathcal{V}_1)\subset\mathcal{V}^{1,1}$ 
is an isotropic vector subspace. Hence $\dim\mathcal{V}_1\leq 1$.  
This gives a contradiction since we would have that $\dim\U(Y_1)\leq 3$.  

Denote by $N_1=\ker\phi_{Y_1}|_{\ker\phi_X}$ and $N_2=\ker\phi_{Y_2}|_{N_1}$.
Then $N_2=\mathcal{N}(\phi)$ since by \eqref{secondst},
\eqref{cod6decomp} and flatness of $\phi$ we have
$$
\lp\phi(Z,\eta),\hat\phi(Y_i,\delta)\rp
=\lp\phi(Z,\eta),\phi(Y_i,\delta)\rp
=\lp\phi(Z,\delta),\phi(Y_i,\eta)\rp=0
$$
for $Z\in\mathfrak{X}(M)$, $\eta\in\Gamma(N_2)$ and
any $\delta\in\Gamma(f_*TM\oplus L_1)$.  It follows from
\eqref{phikerest2} that
$$
\nu(\phi)=\dim N_2\geq \dim N_1-2\geq\dim \ker\phi_X-4
\geq 2n-2p+4\ell_0+2\ell_1,
$$
which is \eqref{phiest1}.\vspace{2ex}\qed

Since $\mathcal{N}(\phi)\cap f_*TM=f_*D_0$ holds, we define 
the $\hat\J_1$-invariant vector subbundle  
$\hat\Lambda\subset f_*D_0^\perp\oplus L_1$  
by the orthogonal decomposition 
\be\label{Lambda2}
\mathcal{N}(\phi)=f_*D_0\oplus\hat\Lambda.
\ee

\begin{proposition}\po\label{charphi} We have that
\be\label{inequalities}
\dim\Sal(\a|_{D\times D})\leq\rank\hat\Lambda\leq\rank L.
\ee 
Moreover, if we have $L=L_0$ and $\rank\hat\Lambda=2\ell$ then $f$ 
has the Kaehler extension 
$F\colon N^{2n+2\ell}\to\R^{2n+p}$ satisfying 
$\nu^c_F=2\ell+d$ given by \eqref{theextension}. 
Furthermore, the extension $F$ is minimal if and only if 
$f$ is minimal.
\end{proposition}

\proof We obtain from \eqref{codhojapluri} that
$$
(\tilde\nabla_X(f_*(\nabla_ST)_{D^\perp}+\a(S,T)))_{L^\perp}
=\a_{L^\perp}(X,\nabla_ST)+(\nabla_X^\perp\a(S,T))_{L^\perp}=0
$$
for any $S,T\in\Gamma(D)$. In particular, since 
$L\subset L_0$ we have 
$$
(\tilde\nabla_X(f_*(\nabla_ST)_{D^\perp}
+\a(S,T)))_{L_0^\perp}=0
$$
for any $X\in\mathfrak{X}(M)$ and $S,T\in\Gamma(D)$.  
Since $D$ is $J$-invariant, then it follows that
$\phi(X,f_*(\nabla_ST)_{D^\perp}+\a(S,T))=0$
for any $X\in\mathfrak{X}(M)$ and $S,T\in\Gamma(D)$. 
Thus 
$$
f_*(\nabla_ST)_{D^\perp}+\a(S,T)\in\Gamma(\hat\Lambda)\;\;
\text{for any}\;\;S,T\in\Gamma(D). 
$$ 
 
Let $\psi\colon N_1^g\to\Sal(\a|_{D\times D})$ the 
isomorphism given by \eqref{defder} and let 
${\cal L}\colon N_1^g\to\hat\Lambda$ be the linear 
map defined by 
${\cal L}(\a^g(S,T))=f_*(\nabla_Si_*T)_{D^\perp}+\a(i_*S,i_*T)$
for any $S,T\in\mathfrak{X}(\Sigma)$. Then \eqref{secondfundg} 
shows that the linear map 
$\mathcal{L}\circ\psi^{-1}\colon\Sal(\a|_{D\times D})\to\hat\Lambda$ 
is injective, thus giving the first inequality in 
\eqref{inequalities}.

We show that $\pi(\mathcal{N}(\phi))\subset L$ where 
$\pi\colon f_*TM\oplus L_1\to L_1$ is the projection
onto the second component. In fact, if
$\xi\in\pi(\mathcal{N}(\phi))$ there is 
$Z\in\mathfrak{X}(M)$ such that $f_*Z+\xi\in\mathcal{N}(\phi)$. 
In particular, we have that
$$
0=\phi(S,f_*Z+\xi)
=((\nabla_S^\perp\xi)_{L_0^\perp},(\nabla_S^\perp\J_1\xi)_{L_0^\perp})
\;\,\text{for any}\;\,S\in\Gamma(D_0).
$$
Then from \eqref{defl2} we obtain that $\xi\in\Gamma(L)$. 
Since $\pi|_{\hat\Lambda}\colon\hat\Lambda\to L$ is injective 
hence also the second inequality of \eqref{inequalities} holds.

For the proof of the remaining part of the statement
we have $L=L_0$. Then
$$
((\tilde\nabla_{X}\delta)_{L^\perp},
(\tilde\nabla_{X}\J_1\delta)_{L^\perp})
=\phi(X,\delta)=0
$$
for any $X\in\mathfrak{X}(M)$ and $\delta\in\mathcal{N}(\phi)$.
Since now $\hat\J=\hat\J_1$, we have that Theorem \ref{develop} 
applies to $\hat\Lambda$ and gives that $f$ has a Kaehler 
extension as stated.\qed

\begin{theorem}\po
Let $f\colon M^{2n}\to\R^{2n+p}$, $p\leq n-1$ and $p=4,6,8,10$,
be a real Kaehler submanifold with  $\nu_f^c(x)<2n-2p$ at any
$x\in M^{2n}$. If $p\geq 8$ assume that $\nu_s^c(x)<2n-2s$ 
at any $x\in M^{2n}$ where $s=6$ if $p=8$ and $s=6,8$ if $p=10$. 
Then $f$ restricted to any connected component of an open dense 
subset of $M^{2n}$ has a Kaehler extension 
$F\colon N^{2n+2r}\to\R^{2n+p}$ satisfying 
$\nu^c_F(z)\geq 2n-2p+6r$ at any $z\in N^{2n+2r}$, where we have 
$(p,r)=(4,1),\,(6,1\,\text{or }2),\, 
(8,2\,\text{or }3)\;\text{or}\;(10,3\,\text{or }4)$.
Moreover, the extension $F$ is minimal if and only if $f$ 
is minimal.
\end{theorem}

\proof We argue for $p=10$ being the other cases 
similar, albeit simpler. Proposition \ref{step3}
gives that $\ell\geq 1$. It also follows from Proposition 
\ref{step3} using \eqref{Gamma2} that $d\geq 2n-20+4\ell$. 
Since $\nu^c_{10-2\ell}\geq\nu^c(\a_{L^\perp})=d$, 
this is for $\ell=1,2$ a contradiction. Hence, we have 
$\ell\geq 3$.

If $\ell_0=5$, then $D_0=TM$ and thus part $(ii)$ of 
Lemma \ref{charK} gives that $f$ is holomorphic, a 
contradiction. Hence $\ell_0\leq 4$ and $\ell_0=\ell_1$
by Proposition \ref{step2}. 
Thus $3\leq\ell\leq\ell_1=\ell_0\leq 4$.

\begin{fact}\po\label{Fact1}  {\em Assume that $\ell_0\geq 3$ 
and that $d_0=2n-20+4\ell_0$. Then $2\ell=2\ell_0$ and $f$ 
has a Kaehler extension $F\colon N^{2n+2\ell}\to\R^{2n+10}$ 
with $\nu^c_F(z)=2n-20+6\ell$ for any $z\in N^{2n+2\ell}$.
}\end{fact}

Proposition \ref{phiestimates} yields \eqref{phiest1}, that is, 
we have $\nu(\phi)\geq 2n-20+6\ell_0$. From \eqref{Lambda2} we 
obtain $\rank\hat\Lambda\geq 2\ell_0$ and then \eqref{inequalities} 
gives that $\rank\hat\Lambda=2\ell=2\ell_0$. Fact \ref{Fact1} 
now follows from Proposition \ref{charphi}.
\vspace{1ex}

Next we analyze the two remaining cases.
\vspace{1ex}

\noindent\emph{Case $\ell_0=3$}. We have that $\ell=\ell_0$ 
and then \eqref{Gamma2} gives $d\geq 2n-8$. If $d=2n-8$ 
we obtain from Fact \ref{Fact1} that $f$ has a Kaehler 
extension $F\colon N^{2n+6}\to\R^{2n+10}$ with $\nu^c_F=2n-2$.

Assume now that $d\geq 2n-6$. We have that 
$\mathcal{S}=\mathcal{S}(\a|_{D\times D})$ satisfies 
$\mathcal{S}=L$ since, if otherwise, then \eqref{c:snullity} 
gives a contradiction with the assumption that $\nu^c_6<2n-12$. 
Now Proposition \ref{lambda} and 
Theorem \ref{develop} yield that $f$ has a Kaehler extension 
$F\colon N^{2n+6}\to\R^{2n+10}$ with $\nu^c_F\geq 2n$.
\vspace{1ex}

\noindent\emph{Case $\ell_0=4$}. From \eqref{Gamma2} we 
have $d_0\geq 2n-4$. If $d_0=2n-4$, then Fact \ref{Fact1}
yields that $f$ has a Kaehler extension
$F\colon N^{2n+8}\to\R^{2n+10}$ with $\nu^c_F=2n+4$.

Assume that $d_0\geq 2n-2$. If $\ell=3$, then Proposition 
\ref{destimate} gives $d\geq 2n-4$. We cannot be in the
situation of part $(i)$ of Lemma \ref{cases} since then
\eqref{c:snullity} yields $\nu^c_6\geq 2n-8$ in contradiction 
with the assumption. Hence, we have that 
$\mathcal{S}=L$, and as above it follows that $f$ has 
a Kaehler extension $F\colon N^{2n+6}\to\R^{2n+10}$ 
with $\nu^c_F\geq 2n+2$.

If $\ell=4$ we cannot have that $d=2n$ since then $N_1^f=L$ 
and condition $(\mathcal{C}_2)$ would give that $f$ is not 
substantial, which has been ruled out. 
Then  by the Fact \eqref{Fact1} we can assume that $d=2n-2$. 
If $L=\Sal$ we have from Proposition \ref{lambda} 
and Theorem \ref{develop} that $f$ has a Kaehler extension 
$F\colon N^{2n+8}\to\R^{2n+10}$ with $\nu^c_F\geq 2n+6$.

Assume that $L\neq\Sal$. Lemma \ref{intkernel} applies to any 
$\J$-invariant vector subbundle 
$\mathcal{L}^{2\bar{r}}\subset\mathcal{R}$ 
and $E=D$ since $A_\eta D\subset D^\perp$ is satisfied for 
any $\eta\in\mathcal{R}$. Then
\be\label{snullest}
\nu^c(\a_{\mathcal{L}\oplus L^\perp})\geq (2n-2)(\bar{r}+1)-2n\bar{r}.
\ee 
It cannot be that $\bar{r}=2$ since otherwise we have
$\rank\mathcal{L}\oplus L^\perp=6$, and then \eqref{snullest} 
gives $\nu^c_6\geq 2n-6$, which is a contradiction.  
Hence $\rank\mathcal{R}=2$, and since
$\mathcal{S}^\perp=\mathcal{R}\oplus L^\perp$ then \eqref{snullest}
yields for $\bar{D}=\mathcal{N}_c(\a_{\Sal^\perp})\subset D$ that
$\rank \bar{D}\geq 2n-4$. 

We claim that the vector subbundle $\Sal$ is parallel along $\bar{D}$
in the normal connection. The  Codazzi equation gives
$$
A_{\psi(T,\eta)}S-A_{\psi(S,\eta)}T=A_\eta[S,T]\in\Gamma(D^\perp),
$$  
where $\psi\colon\Gamma(\bar{D})\times\Gamma(\mathcal{R})\to\Gamma(\Sal)$ 
is defined by $\psi(T,\eta)=(\nabla_T^\perp\eta)_{\Sal}$. 
If $\Sal$ is not parallel along $\bar{D}$, and since
$L$ is parallel along $D$, then there is 
$\eta\in\Gamma(\mathcal{R})$ such that 
$\psi_\eta=\psi(\cdot,\eta)\neq 0$. Let $T\in \bar{D}$ satisfy 
$\psi_\eta T\neq 0$ and $S\in\ker\psi_\eta$. Then 
$$
A_{\psi(T,\eta)}S=A_\eta[S,T]\in\Gamma(D^\perp).
$$
Hence $A_{\Ima\psi_\eta}\ker\psi_\eta\in \Gamma(D^\perp)$. 
By \eqref{cond} then
$A_{\Ima\psi_\eta+\J\Ima\psi_\eta}(\ker\psi_\eta+J\ker\psi_\eta)
\in\Gamma(D^\perp)$. 
Let $\T^2\subset\Ima\psi_\eta+\J\Ima\psi_\eta$ be 
a $\J$-invariant vector subbundle of rank $2$.  Then
Lemma~\ref{intkernel} applied to $\mathcal{L}=\T^2$ and
$E=\ker\psi_\eta+J\ker\psi_\eta$ gives
$$
\nu^c(\a_{\T\oplus\mathcal{R}\oplus L^\perp})
\geq \dim(\ker\psi_\eta+J\ker\psi_\eta)-2.
$$
Since $\dim\ker\psi_\eta\geq\rank \bar{D}-\rank\Sal=2n-10$ then 
$\nu^c_6\geq\nu^c(\a_{\T\oplus\mathcal{R}\oplus L^\perp})\geq 2n-12$,
and this is a contradiction which proves the claim.

Since $\Sal$ is $\J$-invariant it satisfies the condition
$(\mathcal{C}_1)$ and that $\J|_{\Sal}$ is parallel. By the 
claim $\Sal$ also satisfies condition $(\mathcal{C}_2)$.
Moreover, Proposition \ref{lambda} and Theorem \ref{develop}
apply since $(\Sal,\J|_\Sal)$ satisfies the conditions
$(\mathcal{C}_1)$ and $(\mathcal{C}_2)$ in addition to  
$\Sal(\a|_{D_2\times D_2})=\Sal$. 
In fact, if otherwise there is an orthogonal decomposition 
$\Sal=\Sal_1\oplus\mathcal{R}_1$ with $\mathcal{R}_1\neq0$ 
and then \eqref{c:snullity} gives $\nu^c_6\geq 2n-8$.  
It follows that $f$ has a Kaehler extension 
$F\colon N^{2n+6}\to\R^{2n+10}$ with $\nu^c_F\geq 2n+2$.

Finally, the statement about minimality in the theorem holds
since all Kaehler extensions have been obtained by the use 
of either Theorem \ref{develop} or Proposition \ref{charphi}.
\qed

\begin{theorem}\po
Let $f\colon M^{2n}\to\R^{2n+p}$, $p\leq n-1$ and  $p=3,5,7,9,11$, 
be a real Kaehler submanifold with $\nu_f^c(x)<2n-2p$ at any 
$x\in M^{2n}$.  If $p\geq 7$ assume that $\nu_s^c(x)<2n-2s$ at any 
$x\in M^{2n}$ where $s=5$ if $p=7$, $s=5,7$ if $p=9$ and $s=5,7,9$ 
if $p=11$. Then $f$ restricted to any connected component of an 
open dense subset of $M^{2n}$ has a 
Kaehler extension $F\colon N^{2n+2r}\to\R^{2n+p}$ with 
$\nu^c_F(z)\geq 2n-2p+6r$ at any $z\in N^{2n+2r}$, where we have
$(p,r)=(3,1),\, (5,1\,\text{or }2),\,(7,2\,\text{or }3),\, 
(9,3\,\text{or }4)\;\text{or}\;(11,4\,\text{or }5)$.
Moreover, the extension $F$ is minimal if and only if $f$ is minimal.
\end{theorem}

\proof  We argue for $p=11$ being the other cases
similar, albeit simpler. Proposition~\ref{step3}
gives that $\ell\geq 1$. It follows from \eqref{Gamma2} and 
Proposition \ref{step3} that $d\geq 2n-22+4\ell$. Since 
$\nu^c_{11-2\ell}\geq d$, this gives contradiction for $\ell=1,2,3$. 
Hence, we have $\ell\geq 4$. Then we obtain from Proposition 
\ref{step2} that $\ell_0=\ell_1$, and hence 
$4\leq\ell\leq\ell_1=\ell_0\leq 5$.
\vspace{1ex}

\noindent\emph{Case $\ell_0=4$}. By \eqref{Gamma2} we have
$d=d_0\geq 2n-6$. The next Fact, whose proof is roughly the 
same than the one of Fact \ref{Fact1}, gives for $d_0=2n-6$ 
that $f$ has a Kaehler extension 
$F\colon N^{2n+8}\to\R^{2n+11}$ with $\nu^c_F=2n+2$.

\begin{fact}\po\label{Fact2} {\em Assume that $\ell_0\geq 3$ and
that $d_0=2n-22+4\ell_0$. Then $\ell=\ell_0$ and $f$ has a 
Kaehler extension $F\colon N^{2n+2\ell}\to\R^{2n+11}$ with
$\nu^c_F(z)=2n-22+6\ell$ for any $z\in N^{2n+2\ell}$.
}\end{fact}

Assume that $d_0\geq 2n-4$. We have that $L=\Sal=\Sal(\a|_{D\times D})$. 
If otherwise, it follows that $\mathcal{R}\neq 0$ in the 
orthogonal decomposition $L=\Sal\oplus \mathcal{R}$. 
Then $A_\xi D, A_{\J\xi}D\subset D^\perp$ for any
$\xi\in\Gamma(\mathcal{R})$. Since $\ell=\ell_0=4$, we obtain 
from Lemma \ref{intkernel} applied to $\{{\cal R},D\}$ that 
$\nu^c_5\geq 2n-8$, and this is a contradiction. Then 
Theorem \ref{develop} gives a Kaehler extension 
$F\colon N^{2n+8}\to\R^{2n+11}$ with $\nu^c_F\geq 2n+4$.
\vspace{1ex}

\noindent\emph{Case $\ell_0=5$}.  
We can assume that $d_0=2n$ since, if otherwise, 
then \eqref{Gamma2} yields $d_0=2n-2$, and we have 
by Fact \ref{Fact2} a Kaehler extension 
$F\colon N^{2n+10}\to\R^{2n+11}$ with $\nu^c_F=2n+8$.

If $\ell=\ell_0$ then $N_1^f= L$ and we would conclude 
that $f$ is not substantial. Hence $\ell=4$, $\ell_0=5$ 
and $d_0=2n$. Then Proposition \ref{destimate} yields 
$d=2n-2$. It holds that $L=\Sal$ since if $\mathcal{R}\neq 0$ 
then $A_\xi D, A_{\J\xi} D\subset D^\perp$ for
$\xi\in\Gamma(\mathcal{R})$. Hence, we have from Lemma 
\ref{intkernel} that $\nu^c_5\geq 2n-4$, which is
a contradiction. Then Theorem \ref{develop} gives a Kaehler 
extension $F\colon N^{2n+8}\to\R^{2n+11}$ with
$\nu^c_F\geq 2n+6$.

Finally, the statement about minimality in the theorem holds
since all Kaehler extensions have been obtained by the use 
of either Theorem \ref{develop} or Proposition \ref{charphi}.
\qed

\begin{theorem}\po
Let $f\colon M^{2n}\to\R^{2n+p}$, $p\leq n-1$ and $p=7,8$ or 
$9$, be a real Kaehler submanifold with $\nu_f^c(x)<2n-2p$ 
at any $x\in M^{2n}$. Assume that $f$ when restricted to any 
open subset of $M^{2n}$ does not admit Kaehler extension
given \eqref{theextension}. Then $f$ restricted to any 
connected components of an open dense subset of 
$M^{2n}$ is $(2n-t)$-complex ruled for
$(p,t)=(7,8),(8,12)$ or $(9,16)$.  
\end{theorem}

\proof We argue for $p=9$ being the other cases similar, 
albeit simpler. Proposition \ref{step3} then gives
that $1\leq\ell\leq\ell_0\leq 4$.
\vspace{1ex}

\noindent\emph{Case $\ell_0=4$}.
From \eqref{Gamma2} either $d_0=2n$ or $d_0=2n-2$, 
and Proposition \ref{phiestimates} gives that 
$\nu(\phi)\geq 2n-2+2\ell_1$. On the other hand, since 
\eqref{Lambda2} yields  $\nu(\phi)=d_0+\rank\hat\Lambda$, 
then $\rank\hat\Lambda\geq 2\ell_1-2$
if $d_0=2n$ and  $\rank\hat\Lambda\geq 2\ell_1$ if $d_0=2n-2$.

We claim that $d_0=2n-2$. Assume otherwise, that is, that 
$D_0=TM$. By part $(ii)$ of Lemma \ref{charK} we have 
$\ell_1=\ell_0=4$.  Hence $6\leq\rank\hat\Lambda\leq 2\ell$, 
where the second inequality follows from \eqref{inequalities}. 
Then $\ell=3$.  In fact, if otherwise, then $\ell=\ell_0=4$,  
and since $D=TM$, then $N_1^f=L$. Condition $(\mathcal{C}_2)$ 
gives that $f$ reduces codimension, which has been ruled out.  
Now Proposition \ref{destimate} yields that $d\geq 2n-2$, 
which is in contradiction with any of the estimates obtained 
in part $(I\!I\!I)$ of Theorem \ref{ell}, and the claim has 
been proved. 

Since $d_0=2n-2$ from the claim, then 
$2\ell\geq\rank\hat\Lambda\geq2\ell_1$. Hence 
$2\ell=\rank\hat\Lambda=2\ell_1$. Since $f$ does not 
admit a Kaehler extension it follows from Proposition 
\ref{charphi} that $\ell<\ell_0$. Using Proposition 
\ref{step2}, we obtain that $\ell\leq 2$ and 
$d\geq 2n-10+2\ell$. Then either $\ell=1$ and 
$d\geq 2n-8$ or $\ell=2$ and $d\geq 2n-6$. Part $(I)$ of 
Theorem \ref{ell} yields a contradiction for the former 
case. In the latter case, we have from part $(I\!I)$ of 
Theorem \ref{ell} that $f$ is $(2n-12)$-complex ruled.
\vspace{1ex}

\noindent\emph{Case $\ell_0=3$}.
If $\ell_1\leq 2$, then \eqref{Gamma2} and Proposition \ref{step2} 
yield $\ell=\ell_1=1$ and $d\geq 2n-10$. Now part $(I)$ of 
Theorem \ref{ell} gives that $f$ is $(2n-10)$-complex ruled. 

Assume that $\ell_0=\ell_1=3$.
We claim that $d_0\geq 2n-4$.  If otherwise \eqref{Gamma2} gives 
$d_0=2n-6$, and hence Proposition \ref{phiestimates} that 
$\nu(\phi)\geq 2n$. Then \eqref{Lambda2} and \eqref{inequalities} 
yield $6\leq\rank\hat\Lambda\leq2\ell$. Then 
$2\ell=\rank\hat\Lambda=2\ell_0$, and Proposition \ref{charphi} 
gives a contradiction, which proves the claim. 

Suppose that $d_0=2n-4$. Proposition \ref{phiestimates} gives that 
$\nu(\phi)\geq 2n$. Then \eqref{Lambda2} and \eqref{inequalities} 
yield   $4\leq\rank\hat\Lambda\leq2\ell$.  Hence, either $\ell=2$ or 
$\ell=3$. In the former case, Proposition~\ref{destimate} gives 
that $d\geq 2n-6$. Then it follows that items $(b)$ and $(c)$ of 
part $(I\!I)$ of Theorem~\ref{ell} are possible, that is, that
$f$ is $(2n-12)$-complex ruled.
In the latter case, we have $d=2n-4$. Then $f$ is either as
in $(c)$ or $(d)$ or $(e)$ of part 
$(I\!I\!I)$ of Theorem \ref{ell}. We conclude that $f$ 
is $(2n-16)$-complex ruled.

It remains to analyze the case when $d_0\geq 2n-2$. From Proposition 
\ref{destimate} we obtain $d\geq 2n-8+2\ell$. Then either $\ell=1$ 
and $d\geq 2n-6$ or $\ell=2$ and $d\geq 2n-4$ or $\ell=3$ 
and $d\geq 2n-2$, and we have a contradiction with all parts
of Theorem \ref{ell}.
\vspace{1ex}

\noindent\emph{Case $\ell_0=2$}.
We argue that $d_0\geq 2n-8$. If otherwise \eqref{Gamma2} yields
$d_0=2n-10$. By Proposition \ref{step2} we have that $\ell_0=\ell_1$.  
Then Proposition \ref{phiestimates} gives $\nu(\phi)\geq 2n-6$.  
From \eqref{Lambda2} and \eqref{inequalities} we obtain
$4\leq\rank\hat\Lambda\leq2\ell$. Hence 
$\ell=\ell_0=\rank\hat\Lambda$, and Proposition \ref{charphi} 
gives that $f$ has a Kaehler extension, which is a contradiction.

We have that either $\ell=1$ or $\ell=2$. In the former case, 
we have from Proposition~\ref{step2} that $\ell_0=\ell_1$.
Proposition \ref{destimate} gives $d\geq 2n-10$ and part $(I)$ 
of Theorem \ref{ell} that $f$ is $(2n-10)$-complex ruled.  
In the latter case, part $(I\!I)$ of Theorem \ref{ell} 
yields that $f$ is $(2n-16)$-complex ruled.
\vspace{1ex}

\noindent\emph{Case $\ell_0=1$}. We have $\ell=\ell_0$ and 
from \eqref{Gamma2} that $d\geq 2n-14$. Part $(I)$ of 
Theorem \ref{ell} gives that $f$ is $(2n-14)$-complex ruled.
\vspace{4ex}\qed

Marcos Dajczer is partially supported 
by the grant PID2021-124157NB-I00 funded by 
MCIN/AEI/10.13039/501100011033/ `ERDF A way of making Europe',
Spain, and are also supported by Comunidad Aut\'{o}noma de la Regi\'{o}n
de Murcia, Spain, within the framework of the Regional Programme
in Promotion of the Scientific and Technical Research (Action Plan 2022),
by Fundaci\'{o}n S\'{e}neca, Regional Agency of Science and Technology,
REF, 21899/PI/22.

\noindent Sergio J. Chion Aguirre\\
CENTRUM Cat\'{o}lica Graduate Business School,\\
Pontificia Universidad Cat\'{o}lica del Per\'{u},\\
Lima -- Per\'{u},\\
e-mail: sjchiona@pucp.edu.pe
\bigskip

\noindent Marcos Dajczer\\
Departamento de Matem\'{a}ticas\\
Universidad de Murcia,\\
E-30100 Espinardo, Murcia -- Spain\\
e-mail: marcos@impa.br
\end{document}